# The Skorokhod embedding problem and its offspring*

**Jan Obłój**[†]

*Laboratoire de Probabilités et Modèles Aléatoires, Université Paris 6
4 pl. Jussieu, Boîte 188, 75252 Paris Cedex 05, France.
Department of Mathematics, Warsaw University
ul. Banacha 2, 02-097 Warszawa, Poland.
e-mail:* `obloj@mimuw.edu.pl`

**Abstract:** This is a survey about the Skorokhod embedding problem. It presents all known solutions together with their properties and some applications. Some of the solutions are just described, while others are studied in detail and their proofs are presented. A certain unification of proofs, based on one-dimensional potential theory, is made. Some new facts which appeared in a natural way when different solutions were cross-examined, are reported. Azéma and Yor's and Root's solutions are studied extensively. A possible use of the latter is suggested together with a conjecture.

---

*This is an original survey paper.
[†]Work partially supported by French government (scholarship N° 20022708) and Polish Academy of Science.

**Contents**

## 1. Introduction

The so called *Skorokhod embedding problem* or *Skorokhod stopping problem* was first formulated and solved by Skorokhod in 1961 [114] (English translation in 1965 [115]). For a given centered probability measure $\mu$ with finite second moment and a Brownian motion $B$, one looks for an integrable stopping time $T$ such that the distribution of $B_T$ is $\mu$. This original formulation has been changed, generalized or narrowed a great number of times. The problem has stimulated research in probability theory for over 40 years now. A full account of it is hard to imagine; still, we believe there is a point in trying to gather various solutions and applications in one place, to discuss and compare them and to see where else we can go with the problem. This was the basic goal set for this survey and, no matter how far from it we have ended up, we hope that it will prove of some help to other researchers in the domain.

Let us start with some history of the problem. Skorokhod's solution required a randomization external to $B$. Three years later another solution was proposed by Dubins [31] (see Section 3.2), which did not require any external randomization. Around the same time, a third solution was proposed by Root. It was part of his Ph.D. thesis and was then published in an article [103]. This solution has some special minimal properties, and will be discussed in Section 7.

Soon after, another doctoral dissertation was written on the subject by Monroe who developed a new approach using additive functionals (see Section 3.6). His results were published in 1972 [79]. Although he did not have any explicit formulae for the stopping times, his ideas proved fruitful as can be seen from the elegant solutions by Vallois (1983) (see Section 3.12) and Bertoin and Le Jan (1992) (see Section 3.15), which also use additive functionals.

The next landmark was set in 1971 with the work of Rost [105]. He generalized the problem by looking at any Markov process and not just Brownian motion. He gave an elegant necessary and sufficient condition for the existence of a solution to the embedding problem (see Section 3.5). Rost made extensive use of potential theory. This approach was also used a few years later by Chacon and Walsh [19], who proposed a new solution to the original problem which included Dubins' solution as a special case.

By that time, Skorokhod embedding, also called the Skorokhod representation, had been successfully used to prove various invariance principles for random walks (Sawyer [110]). It was the basic tool used to realize a discrete process as a stopped Brownian motion. This ended however with the Komlós, Major and Tusnády [25] construction, which proved far better for this purpose (see Section 11.2). Still, the Skorokhod embedding problem continued to inspire researchers and found numerous new applications (see Sections 10 and 11).

The next development of the theory came in 1979 with a solution proposed by Azéma and Yor [3]. Unlike Rost, they made use of martingale theory, rather than Markov and potential theory, and their solution was formulated for any recurrent, real-valued diffusion. We will see in Section 5 that their solution can be obtained as a limit case of Chacon and Walsh's solution. Azéma and Yor's solution has interesting properties which are discussed as well. In particular the solution maximizes stochastically the law of the supremum up to the stopping time. This direction was continued by Perkins [90], who proposed his own solution in 1985, which in turn minimizes the law of the supremum and maximizes the law of the infimum.

Finally, yet another solution was proposed in 1983 by Bass [4]. He used the stochastic calculus apparatus and defined the stopping time through a time-change procedure (see Section 3.11). This solution is also reported in a recent book by Stroock ([117], p. 213–217).

Further developments can be classified broadly into two categories: works trying to extend older results or develop new solutions, and works investigating properties of particular solutions. The former category is well represented by papers following Monroe's approach: solution with local times by Vallois [118] and the paper by Bertoin and Le Jan [6], where they develop explicit formulae for a wide class of Markov processes. Azéma and Yor's construction was taken as a starting point by Grandits and Falkner [46] and then by Pedersen and Peskir [89] who worked with non-singular diffusions. Roynette, Vallois and Yor [108] used Rost criteria to investigate Brownian motion and its local time [108]. There were also some older works on $n$-dimensional Brownian motion (Falkner [36]) and right processes (Falkner and Fitzsimmons [37]).

The number of works in the second category is greater and we will not try to describe it now, but this will be done progressively throughout this survey. We want to mention, however, that the emphasis was placed on the one hand on the solution of Azéma and Yor and its accurate description and, on the other hand, following Perkins' work, on the control of the maximum and the minimum of the stopped martingale (Kertz and Rösler [61], Hobson [54], Cox and Hobson [23]).

All of the solutions mentioned above, along with some others, are discussed in Section 3, some of them in detail, some very briefly. The exceptions are Azéma and Yor's, Perkins' and Root's solutions, which are discussed in Sections 5, 6 and 7 respectively. Azéma and Yor's solution is presented in detail with an emphasis on a one-dimensional potential theoretic approach, following Chacon and Walsh (cf. Section 3.8). This perspective is introduced, for measures with finite support, in Section 4, where Azéma and Yor's embeddings for random walks are also considered. Perkins' solution is compared with Jacka's solution and both are also discussed using the one-dimensional potential theoretic approach. Section 7, devoted to Root's solution, contains a detailed description of works of Root and Loynes along with results of Rost. The two sections that follow, Sections 8 and 9, focus on certain extensions of the classical case. The former treats the case of non-centered measures and uniform integrability questions and the latter extends the problem and its solutions to the case of real-valued diffusions. The last two sections are concerned with applications: Section 10 develops a link between the Skorokhod embedding problem and optimal stopping theory, and the last section shows how the embedding extends to processes; it also discusses some other applications.

To end this introduction, we need to make one remark about the terminol-ogy. Skorokhod's name is associated with several "problems". In this work, we are concerned solely with the Skorokhod embedding problem, which we call also "Skorokhod's stopping problem" or simply "Skorokhod's problem". As mentioned above, this problem was also once called the "Skorokhod representation" (see for example Sawyer [110], Freedman [42] p. 68) but this name was later abandoned.

## 2. The problem and the methodology

As noted in the Introduction, the Skorokhod embedding problem has been generalized in various manners. Still, even if we stick to the original formulation, it has seen many different solutions, developed using almost all the major techniques which have flourished in the theory of stochastic processes during the last forty years. It is this methodological richness which we want to stress right from the beginning. It has at least two sources. Firstly, different applications have motivated, in the most natural way, the development of solutions formulated within different frameworks and proved in different manners. Secondly, the same solutions, such as the Azéma-Yor solution, were understood better as they were studied using various methodologies.

Choosing the most important members from a family, let alone ordering them, is always a risky procedure! Yet, we believe it has some interest and that it is worthwhile to try. We venture to say that the three most important theoretical tools in the history of Skorokhod embedding have been: a) potential theory, b) the theory of excursions and local times, and c) martingale theory. We will briefly try to argue for this statement, by pointing to later sections where particular embeddings which use these tools are discussed.

As we describe in Section 2.2 below, potential theory has two aspects for us: the simple, one-dimensional potential aspect and the more sophisticated aspect of general potential theory for stochastic processes. The former starts from the study of the potential function $U\mu(x) = -\int |x-y|d\mu(y)$, which provides an convenient description of (centered) probability measures on the real line. It was used first by Chacon and Walsh (see Section 3.8) and proved to be a simple yet very powerful tool in the context of the Skorokhod problem. It lies at the heart of the studies performed by Perkins (see Section 6) and more recently by Hobson [54] and Cox and Hobson (see Section 3.20). We show also how it yields the Azéma-Yor solution (see Section 5.1).

General potential theory was used by Rost to establish his results on the existence of a solution to the Skorokhod embedding problem in a very general setup. Actually, Rost's results are better known as results in potential theory itself, under the name "continuous time filling-scheme" (see Section 3.5). They allowed Rost not only to formulate his own solution but also to understand better the solution developed by Root and its optimal properties (see Section 7). In a sense, Root's solution relies on potential theory, even if he doesn't use it explicitly in his work. Recently, Rost's results were used by Roynette, Vallois and Yor [108] to develop an embedding for Brownian motion and its local time at zero (see Section 3.18).

A very fruitful idea in the history of solving the embedding problem is the following: compare the realization of a Brownian trajectory with the realization of a (function of) some well-controlled increasing process. Use the latter to decide when to stop the former. We have two such "well-controlled" increasing processes at hand: the one-sided supremum of Brownian motion and the local time (at zero). And it comes as no surprise that to realize this idea we need to rely either on martingale theory or on the theory of local times and excursions. The one-sided supremum of Brownian motion appears in Azéma and Yor's solution and their original proof in [3] used solely martingale-theory arguments (see Section 5.4). Local times were first applied to the Skorokhod embedding problem by Monroe (see Section 3.6) and then used by Vallois (see Section 3.12). Bertoin and Le Jan (see Section 3.15) used additive functionals in a very general setup. Excursion theory was used by Rogers [101] to give another proof of the Azéma-Yor solution. It is also the main tool in a recent paper by Obłój and Yor [84].

We speak about a solution to the Skorokhod embedding problem when the embedding recipe works for any centered measure $\mu$. Yet, if we look at a given description of a solution, it was almost always either developed to work with discrete measures and then generalized, or on the contrary, it was designed to work with diffuse measures and only later generalized to all measures. This binary classification of methods happens to coincide nicely with a tool-based classification. And so it seems to us that the solutions which start with atomic measures rely on analytic tools which, more or less explicitly, involve potential theory and, vice versa, these tools are well-suited to approximate any measure by atomic measures. In comparison, the solutions which rely on martingale theory or the theory of excursions and local times, are in most cases easily described and obtained for diffuse measures. A very good illustration for our thesis is the Azéma-Yor solution, as it has been proven using all of the above-mentioned tools. And so, when discussing it using martingale arguments as in Section 5.4, or arguments based on excursion theory (as in Rogers [101]), it is most natural to concentrate on diffuse measures, while when relying on one-dimensional potential theory as in Section 5.1, we approximate any measure with a sequence of atomic measures.

### *2.1. Preliminaries*

We gather in this small subsection some crucial results which will be used throughout this survey without further explanation. $B = (B_t : t \geq 0)$ always denotes a one-dimensional Brownian motion, defined on a probability space $(\Omega, \mathcal{F}, \mathbb{P})$ and $\mathcal{F}_t = \sigma(B_s : s \leq t)$ is its natural filtration, taken completed. The one-sided supremum process is $S_t = \sup_{u \leq t} B_u$ and the two-sided supremum process is $B_t^* = \sup_{u \leq t} |B_u|$. Let $T$ denote an arbitrary stopping time, $T_x = \inf\{t \geq 0 : B_t = x\}$ the first hitting time of $x$, and $T_{a,b} = \inf\{t \geq 0 : B_t \notin (a,b)\}$, $a < 0 < b$, the first exit time form $(a,b)$. $M = (M_t : t \geq 0)$ will always be a real-valued, continuous local martingale and $(\langle M \rangle_t : t \geq 0)$ its quadratic variation process. When no confusion is possible, its one-sided and two-sided suprema, and first exit times will be also denoted $S_t$, $M_t^*$ and $T_{a,b}$ respectively.

For a random variable $X$, $\mathcal{L}(X)$ denotes its distribution and $X \sim \mu$, or $\mathcal{L}(X) = \mu$, means "$X$ has the distribution $\mu$." The support of a measure $\mu$ is denoted $supp(\mu)$. For two random variables $X$ and $Y$, $X \sim Y$ and $X \stackrel{\mathcal{L}}{=} Y$ signify both that $X$ and $Y$ have the same distribution. The Dirac point mass at $a$ is denoted $\delta_a$, or $\delta_{\{a\}}$ when the former might cause confusion. The Normal distribution with mean $m$ and variance $\sigma^2$ is denoted $\mathcal{N}(m, \sigma^2)$. Finally, $\mu_n \Rightarrow \mu$ signifies weak convergence of probability measures.

Note that as $(B_{t \wedge T_{a,b}})_{t \geq 0}$ is bounded, it is a uniformly integrable martingale and the optional stopping theorem (see Revuz and Yor [99], pp. 70–72) yields $\mathbb{E} B_{T_{a,b}} = 0$. This readily implies that $\mathbb{P}(B_{T_{a,b}} = a) = \frac{b}{b-a} = 1 - \mathbb{P}(B_{T_{a,b}} = b)$. This generalizes naturally to situations when Brownian motion does not start from zero. This property is shared by all continuous martingales and, although it looks very simple, it is fundamental for numerous works discussed in this survey, as will be seen later (cf. Section 3.8).

The process $B_t^2 - t$ is easily seen to be a martingale. By optional stopping theorem we see that if $\mathbb{E} T < \infty$ then $\mathbb{E} T = \mathbb{E} B_T^2$. We can actually prove a stronger result (see Root [103] and Sawyer [110]):

**Proposition 2.1.** *Let $T$ be a stopping time such that $(B_{t \wedge T})_{t \geq 0}$ is a uniformly integrable martingale. Then there exist universal constants $c_p, C_p$ such that*

$$c_p \mathbb{E}\left[T^{p/2}\right] \leq \mathbb{E}\left[|B_T|^p\right] \leq C_p \mathbb{E}\left[T^{p/2}\right] \quad \text{for } p > 1. \tag{2.1}$$

*Proof.* The Burkholder-Davis-Gundy inequalities (see Revuz and Yor [99] p. 160) guarantee existence of universal constants $k_p$ and $K_p$ such that $k_p \mathbb{E} T^{p/2} \leq \mathbb{E}[(B_T^*)^p] < K_p \mathbb{E} T^{p/2}$, for any $p > 0$. As $(B_{t \wedge T} : t \geq 0)$ is a uniformly integrable martingale we have $\sup_t \mathbb{E} |B_{t \wedge T}|^p = \mathbb{E} |B_T|^p$, and Doob's $L^p$ inequalities yield $\mathbb{E} |B_T|^p \leq \mathbb{E}[(B_T^*)^p] \leq \left(\frac{p}{p-1}\right)^p \mathbb{E} |B_T|^p$ for any $p > 1$. The proof is thus completed taking $c_p = k_p \left(\frac{p}{p-1}\right)^p$ and $C_p = K_p$. □

We note that the above Proposition is not true for $p = 1$. Indeed we can build a stopping time $T$ such that $(B_{t \wedge T} : t \geq 0)$ is a uniformly integrable martingale, $\mathbb{E} |B_T| < \infty$, and $\mathbb{E} \sqrt{T} = \infty$ (see Exercise II.3.15 in Revuz and Yor [99]).

The last tool that we need to recall here is the Dambis-Dubins-Schwarz theorem (see Revuz and Yor [99], p. 181). It allows us to see a continuous local martingale $(M_t : t \geq 0)$ as a time-changed Brownian motion. Moreover, the time-change is given explicitly by the quadratic variation process $(\langle M \rangle_t : t \geq 0)$. This way numerous solutions of the Skorokhod embedding problem for Brownian motion can be transferred to the setup of continuous local martingales.

## 2.2. On potential theory

Potential theory has played a crucial role in a number of works about the Skorokhod embedding problem and we will also make a substantial use of it, especially in its one-dimensional context. We introduce here some potential the-oretic objects which will be important for us. As this is not an introduction to Markovian theory we will omit certain details and assumptions. We refer to the classical work of Blumenthal and Getoor [9] for precise statements.

On a locally compact space $(E, \mathcal{E})$ with a denumerable basis and a Borel $\sigma$-field, consider a Markov process $X = (X_t : t \geq 0)$ associated with $(P_t^X : t \geq 0)$, a standard semigroup of submarkovian kernels. A natural interpretation is that $\nu P_t^X$ represents the law of $X_t$, under the starting distribution $X_0 \sim \nu$. Define the potential kernel $U^X$ through $U^X = \int_0^\infty P_t^X dt$. This can be seen as a linear operator on the space of measures on $E$. The intuitive meaning is that $\nu U^X$ represents the occupation measure for $X$ along its trajectories, where $X_0 \sim \nu$.

If the potential operator is finite[1] it is not hard to see that for two bounded stopping times, $S \leq T$, we have $U^X(P_S^X - P_T^X) \geq 0$. This explains how the potential can be used to keep track of the relative stage of the development of the process (see Chacon [17]). We will continue the discussion of general potential theory in Section 3.5 with the works of Rost.

For $X = B$, a real-valued Brownian motion, as the process is recurrent, the potential $\nu U^B$ is infinite for positive measures $\nu$. However, if the measure $\nu$ is a signed measure with $\nu(\mathbb{R}) = 0$, and $\int |x|\nu(dx) < \infty$, then the potential $\nu U^B$ is not only finite but also absolutely continuous with respect to Lebesgue measure. Its Radon-Nikodym derivative is given by

$$\frac{d(\nu U^B)}{dx}(x) = -\int |x - y|\nu(dy). \tag{2.2}$$

The RHS is well defined for any probability measure $\mu$ on $\mathbb{R}$ (instead of $\nu$) with $\int |x|d\mu(x) < \infty$ and, with a certain abuse of terminology, this quantity is called the one-dimensional potential of the measure $\mu$:

**Definition 2.2.** *Denote by $\mathcal{M}^1$ the set of probability measures on $\mathbb{R}$ with finite first moment, $\mu \in \mathcal{M}$ iff $\int |x|d\mu(x) < \infty$. Let $\mathcal{M}_m^1$ denote the subset of measures with expectation equal to $m$, $\mu \in \mathcal{M}_m^1$ iff $\int |x|d\mu(x) < \infty$ and $\int x d\mu(x) = m$. Naturally $\mathcal{M}^1 = \bigcup_{m \in \mathbb{R}} \mathcal{M}_m^1$. The one-dimensional potential operator $U$ acting from $\mathcal{M}^1$ into the space of continuous, non-positive functions, $U : \mathcal{M}^1 \to C(\mathbb{R}, \mathbb{R}_-)$, is defined through $U\mu(x) := U(\mu)(x) = -\int_{\mathbb{R}} |x-y|\mu(dy)$. We refer to $U\mu$ as to **the potential of** $\mu$.*

We adopt the notation $U\mu$ to differentiate this case from the general case of a potential kernel $U^X$ when $\mu U^X$ is a measure and not a function. This simple operator enjoys some remarkable properties, which will be crucial for the Chacon-Walsh methodology (see Section 3.8), which in turn is our main tool in this survey. The following can be found in Chacon [17], and Chacon and Walsh [19]:

**Proposition 2.3.** *Let $m \in \mathbb{R}$ and $\mu \in \mathcal{M}_m^1$. Then*

 *(i) $U\mu$ is concave and Lipschitz-continuous with parameter $1$;*

---
[1] That is, for any finite starting measure, it provides a finite measure.

(ii) $U\mu(x) \leq U\delta_{\{m\}}(x) = -|x-m|$ and if $\nu \in \mathcal{M}^1$ and $U\nu \leq U\mu$ then $\nu \in \mathcal{M}^1_m$;

(iii) for $\mu_1, \mu_2 \in \mathcal{M}^1_m$, $\lim_{|x|\to\infty} |U\mu_1(x) - U\mu_2(x)| = 0$;

(iv) for $\mu_n \in \mathcal{M}^1_m$, $\mu_n \Rightarrow \mu$ if and only if $U\mu_n(x) \xrightarrow[n\to\infty]{} U\mu(x)$ for all $x \in \mathbb{R}$;

(v) for $\nu \in \mathcal{M}^1_0$, if $\int_\mathbb{R} x^2 \nu(dx) < \infty$, then $\int_\mathbb{R} x^2 \nu(dx) = \int_\mathbb{R} \Big||x| + U\nu(x)\Big| dx$;

(vi) for $\nu \in \mathcal{M}^1_m$, $U\nu|_{[b,\infty)} = U\mu|_{[b,\infty)}$ if and only if $\mu|_{(b,\infty)} \equiv \nu|_{(b,\infty)}$;

(vii) let $B_0 \sim \nu$, and define $\rho$ through $\rho \sim B_{T_{a,b}}$, then $U\rho|_{(-\infty,a]\cup[b,\infty)} = U\nu|_{(-\infty,a]\cup[b,\infty)}$ and $U\rho$ is linear on $[a,b]$.

*Proof.* We only prove *(ii)*, *(vi)* and *(vii)*. Let $\mu \in \mathcal{M}^1_m$. The first assertion in *(ii)* follows from Jensen's inequality as

$$U\mu(x) = -\int_{-\infty}^{\infty} |x-y| d\mu(y) \leq -\Big|\int_{-\infty}^{\infty} (x-y) d\mu(y)\Big| = U\delta_m(x) = -|x-m|.$$

To prove the other assertions we rewrite the potential:

$$\begin{aligned} U\mu(x) &= -\int_\mathbb{R} |x-y| d\mu(y) = -\int_{(-\infty,x)} (x-y) d\mu(y) - \int_{[x,\infty)} (y-x) d\mu(y) \\ &= x\Big(2\mu([x,\infty)) - 1\Big) + m - 2\int_{[x,\infty)} y d\mu(y), \end{aligned} \quad (2.3)$$

where, to obtain the third equality, we use the fact that $\mu \in \mathcal{M}^1_m$ and so $\int_{[x,\infty)} y d\mu(y) = m - \int_{(-\infty,x)} y d\mu(y)$. The second assertion in *(ii)* now follows. For $\nu \in \mathcal{M}^1$ with $U\nu \leq U\mu$ we have $U\nu(x) \leq -|x-m|$. Using the expression (2.3) above and letting $x \to -\infty$ we see that $\int_\mathbb{R} x d\nu(x) \geq m$ and, likewise, letting $x \to \infty$ we see that the reverse holds.

The formula displayed in (2.3) shows that the potential is linear on intervals $[a,b]$ such that $\mu((a,b)) = 0$. Furthermore, it shows that for $\nu \in \mathcal{M}^1_m$, $U\nu|_{[b,\infty)} = U\mu|_{[b,\infty)}$ if and only if $\mu|_{(b,\infty)} \equiv \nu|_{(b,\infty)}$. Note that the same is true with $[b,\infty)$ replaced by $(-\infty,a]$ and that in particular $U\nu \equiv U\mu$ if and only if $\nu \equiv \mu$. The assertions in *(vii)* follow from the fact that if $B_0 \sim \nu$, then the law of $B_{T_{a,b}}$ coincides with $\nu$ on $(-\infty,a)\cup(b,\infty)$ and does not charge the interval $(a,b)$. □

### 2.3. The Problem

Let us state formally, in its classical form, the main problem considered in this survey.

**Problem (The Skorokhod embedding problem).** *For a given probability measure $\mu$ on $\mathbb{R}$, such that $\int_\mathbb{R} |x| d\mu(x) < \infty$ and $\int_\mathbb{R} x d\mu(x) = 0$, find a stopping time $T$ such that $B_T \sim \mu$ and $(B_{t\wedge T} : t \geq 0)$ is a uniformly integrable martingale.*

Actually, in the original work of Skorokhod [115], the variance $v = \int_{\mathbb{R}} x^2 d\mu(x)$, was assumed to be finite, $v < \infty$, and the condition of uniform integrability of $(B_{t \wedge T} : t \geq 0)$ was replaced by $\mathbb{E} T = v$. Naturally the latter implies that $(B_{t \wedge T} : t \geq 0)$ is a uniformly integrable martingale. However the condition $v < \infty$ is somewhat artificial and the condition of uniform integrability of $(B_{t \wedge T} : t \geq 0)$ seems the correct way to express the idea that $T$ should be "small" (see also Section 8). We point out that without any such condition on $T$, there is a trivial solution to the above problem. We believe it was first observed by Doob. For any probability measure $\mu$, we define the distribution function $F_\mu(x) = \mu((-\infty, x])$ and take $F_\mu^{-1}$ its right-continuous inverse. The distribution function of a standard Normal variable is denoted $\Phi$. The stopping time $T = \inf\{t \geq 2 : B_t = F_\mu^{-1}(\Phi(B_1))\}$ then embeds $\mu$ in Brownian motion, $B_T \sim \mu$. However, we always have $\mathbb{E} T = \infty$.

Skorokhod first formulated and solved the problem in order to be able to embed random walks into Brownian motion (see Section 11.1). We mention that a similar question was posed earlier by Harris [52]. He described an embedding of a simple random walk on $\mathbb{Z}$ (inhomogeneous in space) into Brownian motion by means of first exit times. Harris used the embedding to obtain some explicit expressions for mean recurrence and first-passage times for the random walk.

## 3. A guided tour through the different solutions

In Section 2 above, we tried to give some classification of known solutions to the Skorokhod embedding problem based mainly on their methodology. We intended to describe the field and give the reader some basic intuition. We admit that our discussion left out some solutions, as for example the one presented by Bass (see Section 3.11), which did not fit into our line of reasoning. We will find them now, as we turn to the chronological order. We give a brief description of each solution, its properties and relevant references. The references are summarized in Table 1. Stopping times, solutions obtained by various authors, are denoted by $T_{\cdot}$ with subscript representing author's last name, i.e. Skorokhod's solution is denoted by $T_\mathbf{S}$.

### 3.1. Skorokhod (1961)

In his book from 1961 Skorokhod [114] (English translation in 1965 [115]) defined the problem under consideration and proposed the first solution. He introduced it to obtain some invariance principles for random walks. Loosely speaking, the idea behind his solution is to observe that any centered measure with two atoms can be embedded by the means of first exit time. One then describes the measure $\mu$ as a mixture of centered measures with at most two atoms and the stopping rule consists in choosing independently one of these measures, according to the mixture, and then using the appropriate first exit time. Actually Skorokhod worked with measures with density. His ideas were then extended to work with arbitrary centered probability measures by various authors. We mention Strassen [116], Breiman [11] (see Section 3.4) and Sawyer [109]. We present here a rigorous solution found in Freedman's textbook ([42], pp. 68–76) and in Perkins [90]. Freedman used the name "Skorokhod representation" rather than "Skorokhod embedding problem" alluding to the motivation of representing a random walk as a Brownian motion stopped with a family of stopping times. He used this representation to prove Donsker's theorem and Strassen's law (see Section 11.2). He also obtained some technical lemmas, which he exploited in proofs of invariance principles for Markov chains ([43], pp. 83–87).

For Brownian motion $(B_t)$ and a centered probability measure $\mu$ on $\mathbb{R}$, define for any $\lambda > 0$

$$-\rho(\lambda) = \inf\left\{y \in \mathbb{R} : \int_\mathbb{R} \mathbf{1}_{\{x \leq y\} \cup \{x \geq \lambda\}} x d\mu(x) \leq 0\right\}.$$

Let $R$ be an independent random variable with the following distribution function:

$$\mathbb{P}(R \leq x) = \int_\mathbb{R} \mathbf{1}_{y \leq x}\left(1 + \frac{y}{\rho(y)}\right) d\mu(y).$$

Then, the stopping time defined by

$$T_\mathbf{S} = \inf\{t : B_t \notin (-\rho(R), R)\},$$

where **S** stands for "Skorokhod", satisfies $B_{T_\mathbf{S}} \sim \mu$ and $(B_{t \wedge T_\mathbf{S}} : t \geq 0)$ is a uniformly integrable martingale. In his original work Skorokhod assumed moreover that the second moment of $\mu$ was finite but this is not necessary.

Sawyer [109] investigated in detail the properties of constructed stopping times under exponential integrability conditions on the target distribution $\mu$. In particular he showed that if for all $x > 0$, $\mu((-\infty, x] \cup [x, \infty)) \leq C \exp(-ax^\epsilon)$, for some positive $C, a, \epsilon$, then $\mathbb{P}(T_\mathbf{S} \geq x) \leq C_1 \exp(-bx^\delta)$, with $b = a^{1-a}$, $\delta = \epsilon/(2+\epsilon)$ and some positive constant $C_1 > 0$. One has then $\mathbb{E}(\exp(bT_\mathbf{S}^\delta)) \leq C_2 \int_0^\infty \exp(a|x|^\epsilon) d\mu(x)$, with $C_2$ some positive constant depending only on $\epsilon$.

### 3.2. Dubins (1968)

Three years after the translation of Skorokhod's original work, Dubins [31] presented a different approach which did not require any additional random variable independent of the original process itself. The idea of Dubins is crucial as it will be met in various solutions later on. It is clearly described in Meyer [76]. Given a centered probability measure $\mu$, we construct a family of probability measures $\mu_n$ with finite supports, converging to $\mu$. We start with $\mu_0 = \delta_0$, the measure concentrated on the expectation of $\mu$ (i.e. 0). We look then at $\mu_-$ and $\mu_+$ - the measures resulting from restricting (and renormalizing) $\mu$ to two intervals, $(-\infty, 0]$ and $(0, +\infty)$. We take their expectations: $m_+$ and $m_-$. Choosing the weights in a natural way we define a centered probability measure $\mu_1$ (concentrated on 2 points $m_+$ and $m_-$). The three points considered so far $(m_-, 0, m_+)$ cut the real line into 4 intervals and we again look at renormalized restrictions of $\mu$ on these 4 intervals. We take their expectations, choose appropriate weights and obtain a centered probability measure $\mu_2$. We continue this construction: in the $n^{th}$ step we pass from a measure $\mu_n$ concentrated on $2^{n-1}$ points and a total of $2^n - 1$ points cutting the real line, to the measure $\mu_{n+1}$ concentrated on $2^n$ new points (and we have therefore a total of $2^{n+1} - 1$ points cutting the real line). Not only does $\mu_n$ converge to $\mu$ but the construction embeds easily in Brownian motion using the successive hitting times $(T_n^k)_{1 \leq k \leq 2^{n-1}}$ of pairs of new points (nearest neighbors in the support of $\mu_{n+1}$ of a point in the support of $\mu_n$). $T$ is the limit of $T_n^{2^{n-1}}$ and it is straightforward to see that it solves the Skorokhod problem.

Bretagnolle [12] has investigated this construction and in particular the exponential moments of Dubins' stopping time. Dubins' solution is also reported in Billingsley's textbook [7]. However his description seems overly complicated[2].

### 3.3. Root (1968)

The third solution to the problem came from Root [103]. We will devote Section 7 to discuss it. This solution has a remarkable property - it minimizes the variance (or equivalently the second moment) of the stopping time, as shown by Rost [107]. Unfortunately, the solution is not explicit and is actually hard to describe even for simple target distributions.

### 3.4. Hall (1968)

Hall described his solution to the original Skorokhod problem in a technical report at Stanford University [49], which seems to have remained otherwise unpublished. He did however publish an article on applications of the Skorokhod embedding to sequential analysis [50], where he also solved the problem for Brownian motion with drift. We stress this, as the problem was treated again over 30 years later by Falkner, Grandits, Pedersen and Peskir (cf. Section 3.17) and Hall's work is never cited. Hall's solution to the Skorokhod embedding for Brownian motion with a drift was an adaptation of his method for standard Brownian motion, in very much the same spirit as the just mentioned authors were adapting Azéma and Yor's solution somewhat 30 years later.

Hall [49] looked at Skorokhod's and Dubins' solutions and saw that one could start as Dubins does but then proceed in a way similar to Skorokhod, only this time making the randomization internal.

More precisely, first he presents an easy randomized embedding based on an independent, two-dimensional random variable. For $\mu$, an integrable probability measure on $\mathbb{R}$, $\int |x| d\mu(x) < \infty$, $\int x d\mu(x) = m$, define a distribution on $\mathbb{R}^2$ through $dH_\mu(u, v) = (v - u)\left(\int_m^\infty (x - m) d\mu(x)\right)^{-1} \mathbf{1}_{u \leq m \leq v} d\mu(u) d\mu(v)$. Then, for a centered measure ($m = 0$), the stopping time $T = \inf\{t \geq 0 : B_t \notin (U, V)\}$, where $(U, V) \sim H$ independent of $B$, embeds $\mu$ in $B$, i.e. $B_T \sim \mu$. This solution was first described in Breiman [11]. We remark that the idea is similar to the original Skorokhod's embedding and its extension by Freedman and Perkins (cf. Section 3.1), namely one represents a measure $\mu$ as a random mixture of measure with two atoms. The difference is that in Skorokhod's construction $V$ was a deterministic function of $U$.

---

[2]As a referee has pointed out, the justification of the convergence of the construction, on page 517, could be considerably simplified through an argument given in Exercise II.7 in Neveu [81], p. 37.

The second step consists in starting with an exit time, which will give enough randomness in $B$ itself to make the above randomization internal. The crucial tool is the classical Lévy result about the existence of a measurable transformation from a single random variable with a continuous distribution to a random vector with an arbitrary distribution. Define $\mu_+$ and $\mu_-$, as in Section 3.2, to be the renormalized restrictions of $\mu$ to $\mathbb{R}_+$ and $\mathbb{R}_-$ respectively, and denote their expectations $m_+$ and $m_-$ (i.e. $m_+ = \frac{1}{\mu((0,\infty))} \int_{(0,\infty)} x d\mu(x)$, $m_- = \frac{1}{\mu((-\infty,0])} \int_{(-\infty,0]} x d\mu(x)$). Take $T_1 = T_{m_+ \wedge m_-}$ to be the first exit time from $[m_-, m_+]$. It has a continuous distribution. There exist therefore measurable transformations $f^+$ and $f^-$ such that $(U^i, V^i) = f^i(T_1)$ has distribution $H_{\mu_i}$, $i \in \{+, -\}$. Working conditionally on $B_{T_1} = m_i$, the new starting (shifted) Brownian motion $W_s = B_{s+T_1}$, and the stopping time $T_1$ are independent. Therefore, the stopping time $T_2 = \inf\{s \geq 0 : W_s \notin (U^i, V^i)\}$ embeds $\mu_i$ in the new (shifted) Brownian motion $W$. Finally then, the stopping time $T_{\mathbf{H}} = T_1 + T_2$ embeds $\mu$ in $B$.

As stated above, in [50] this solution was adapted to serve in the case of $B_t^{(\delta)} = B_t + \delta t$ and Hall identified correctly the necessary and sufficient condition on the probability measure $\mu$ for the existence of an embedding, namely $\int e^{-2\delta x} \mu(dx) \leq 1$. His method is very similar to later works of Falkner and Grandits (cf. Section 3.17 and Section 9). As a direct application, Hall also described the class of discrete submartingales that can be embedded in $(B_t^{(\delta)} : t \geq 0)$ by means of an increasing family of stopping times.

### 3.5. Rost (1971)

The first work [104], written in German, was generalized in [105]. Rost's main work on the subject [105], does not present a solution to the Skorokhod embedding problem in itself. It is actually best known for its potential theoretic results, known as the "filling scheme" in continuous time. It answers however a general question about the existence of a solution to the embedding problem for an arbitrary Markov process. The formulation in terms of a solution to the Skorokhod embedding problem is found in [106].

We will start with an attempt to present the general continuous-time filling scheme and the criterion on the laws of stopped processes it yields. We will then specialize to the case of transient processes. We will close this section with a remark on a link between Rost's work and some discrete time works (e.g. Pitman [97]). We do not pretend to be exhaustive and sometimes we make only formal calculations, as our basic aim is just to give an intuitive understanding of otherwise difficult and involved results.

Consider, as in Section 2.2, a Markov process $(X_t)_{t\geq 0}$ on $(E, \mathcal{E})$, with potential operator $U^X$. Let $\nu, \mu$ be two positive measures on $\mathcal{E}$, with finite variation. $\nu$ is thought of as the starting measure of $X$, $X_0 \sim \nu$, and $\mu$ is the measure for which we try to determine whether it can be the distribution of $X$ stopped at some stopping time[3]. Rost shows that there exists a sequence of measures $(\nu_t)$, weakly right continuous in $t$, and a stopping time $S$ (possible randomized) such that:

- $0 \leq \nu_{t+s} \leq \nu_t P_s^X \leq \nu P_{t+s}^X, \quad \forall\, t, s > 0$;
- $\nu_t$ is the law of $X_t \mathbf{1}_{S>t}$.

The family $\mu_t$ is defined through

$$\mu - \mu_t + \nu_t = \nu P_{S \wedge t}^X, \quad t \geq 0. \tag{3.1}$$

We set $\mu_\infty = \downarrow \lim_{t \to \infty} \mu_t$, $\tilde{\mu} = \mu - \mu_\infty$. We say that $\nu \prec \mu$ ("$\nu$ is earlier than $\mu$") if and only if $\mu_\infty = 0$. This defines a partial order on positive measures with finite variation on $\mathcal{E}$.

**Theorem 3.1 (Rost [105]).** *We have $\nu \prec \mu$ if and only if $\mu = \nu P_T^X$ for some (possibly randomized) stopping time $T$.*

This gives a complete description of laws that can be embedded in $X$. However this description, as presented here, may be hard to grasp and we specialize to the case of transient Markov processes to explain it. Luckily, one can pass from one case to another. Namely, given a Markov process $X$ we can consider $X^\alpha$, a transient Markov process, given by $X$ killed at an independent time with exponential distribution (with parameter $\alpha$). A result of Rost assures that if $\mu_\infty^\alpha$ is the measure $\mu_\infty$ constructed as above but with $X$ replaced by $X^\alpha$, then $\mu_\infty = \downarrow \lim_{\alpha \to 0} \mu_\infty^\alpha$.

We suppose therefore that $X$ is transient, and that $\mu U^X$ and $\nu U^X$ are $\sigma$-finite. For an almost Borel set $\mathcal{A} \subset E$ the first hitting time of $\mathcal{A}$, $T_\mathcal{A} := \inf\{t \geq 0 : X_t \in \mathcal{A}\}$, is a stopping time. Note that the measure $\nu P_{T_\mathcal{A}}^X$ is just the law of $X_{T_\mathcal{A}}$.

Recall, that the réduite of a measure $\rho$ is the smallest excessive measure which majorizes $\rho$. It exists, is unique, and can be written as $P_{T_\mathcal{A}}^X \nu$ for a certain $\mathcal{A} \subset E$. We have the following characterization of the filling scheme, or balayage order (see Meyer [77], Rost [105]).

*The measure $\mu$ admits a decomposition $\mu = \overline{\mu} + \mu_\infty$, where $\mu_\infty U^X$ is the réduite of $(\mu - \nu) U^X$ and $\overline{\mu}$ is of the form $\nu P_T^X$, $T$ a stopping time. Furthermore, there exists a finely closed set $\mathcal{A}$ which carries $\mu_\infty$ and for which $\mu_\infty = (\mu - \nu) P_{T_\mathcal{A}}^X$ or $(\overline{\mu} - \nu) P_{T_\mathcal{A}}^X = 0$. In the special case $\mu U^X \leq \nu U^X$ we have $\mu_\infty = 0$ and $\mu = \nu P_T^X$ for some $T$.*

We can therefore reformulate Theorem 3.1 above: for some initial distribution $X_0 \sim \nu$, there exists a stopping time $T$ such that $X_T \sim \mu$ if and only if

$$\mu U^X \leq \nu U^X. \tag{3.2}$$

---

[3] Note that in our notation the roles of $\mu$ and $\nu$ are switched, compared to the articles of Rost [105, 107].

We now comment on the above: it can be rephrased as "$\mu$ happens after $\nu$". It says that
$$\mu U^X(f) \leq \nu U^X(f),$$
for any positive, Borel function $f : E \to \mathbb{R}_+$, where $\nu U^X(f) = \mathbb{E}_\nu \left[ \int_0^\infty f(X_t)dt \right]$, and $\mathbb{E}_\nu$ denotes the expectation under the initial measure $\nu$, $X_0 \sim \nu$. The RHS can be written as:
$$\nu U^X(f) = \mathbb{E}_\nu \left[ \int_0^T f(X_t)dt + \int_T^\infty f(X_t)dt \right] = \mathbb{E}_\nu \left[ \int_0^T f(X_t)dt \right] + \mu U^X(f).$$

One sees that in order to be able to define $T$, condition (3.2) is needed. On the other hand if it is satisfied, one can also hope that $\mathbb{E}_\nu \left[ \int_0^T f(X_t)dt \right] = \nu U^X(f) - \mu U^X(f)$, defines a stopping time. Of course this is just an intuitive justification of a profound result and the proof makes a considerable use of potential theory.

Rost developed the subject further and formulated it in the language of the Skorokhod embedding problem in [106]. We note that the stopping rules obtained in this way may be randomized as in the original work of Skorokhod.

An important continuation of Rost's work is found in two papers of Fitzsimmons. In [38] (with some minor correction recorded in [39]) an existence theorem for embedding in a general Markov process is obtained as a by-product of investigation of a natural extention of Hunt's balayage $L_B m$ (see [38] for appropriate definitions). The second paper [40], is devoted to the Skorokhod problem for general right Markov processes. Suppose that $X_0 \sim \nu$ and the potentials of $\nu$ and $\mu$ are $\sigma$-finite and $\mu U^X \leq \nu U^X$. It is shown then that there exists then a monotone family of sets $\{C(r); 0 \leq r \leq 1\}$ such that if $T$ is the first entrance time of $C(R)$, where $R$ is independent of $X$ and uniformly distributed over $[0, 1]$, then $X_T \sim \mu$. Moreover, in the case where both potentials are dominated by some multiple of an excessive measure $m$, then the sets $C(r)$ are given explicitly in terms of the Radon-Nikodym derivatives $\frac{d(\nu U^X)}{dm}$ and $\frac{d(\mu U^X)}{dm}$.

To close this section, we make a link between the discrete filling scheme and some other works for discrete Markov chains. Namely, we consider the occupation formula for Markov chains developed by Pitman ([96], [97]). Let $(X_n : n \geq 0)$ be a homogenous Markov chain on some countable state space $J$, with transition probabilities given by $P = P(i,j)_{i,j \in J}$. Let $T$ be a finite stopping time for $X$. Define the pre-$T$ occupation measure $\gamma_T$ through $\gamma_T(A) = \sum_{n=0}^\infty \mathbb{P}(X_n \in A, T > n)$, $A \subset J$. Let $\lambda$ be the initial distribution of the process, $X_0 \sim \lambda$, and $\lambda_T$ the distribution of $X_T$. Then Pitman [97], as an application of the occupation measure formula, obtains:

$$\gamma_T + \lambda_T - \lambda = \gamma_T P. \tag{3.3}$$

Given a distribution $\mu$ on $J$ this can serve to determine whether there exists a stopping time $T$ such that $\lambda_T = \mu$.

We point out that the above can be obtained via Rost's discrete filling scheme. Rost ([105] p.3) obtains

$$\tilde{\mu} - \lambda + M = MP, \quad \text{where } M = \sum_n \lambda_n. \tag{3.4}$$

Applying this to the process $X^T = (X_{n \wedge T} : n \geq 0)$, we have $M = \gamma_T$. Rost's condition in Theorem 3.1, namely $\tilde{\mu} = \mu$, is then equivalent to (3.3).

The occupation measure formula admits a continuous-time version (see Fitzsimmons and Pitman [41]), but we will not go into details nor into its relations with the continuous-time filling scheme of Rost.

### 3.6. Monroe (1972)

We examine this work only very briefly. Monroe [79] doesn't have any explicit solutions but again what is important is the novelty of his methodology. Let $(L_t^a)_{t \geq 0}$ be the local time at level $a$ of Brownian motion. Then the support of $dL_t^a$ is contained in the set $\{B_t = a\}$ almost surely. We can use this fact and try to define the stopping time $T$ through inequalities involving local times at levels $a$ ranging through the support of the desired terminal distribution $\mu$. The distribution of $B_T$ will then have the same support as $\mu$. Monroe [79] shows that we can carry out this idea and define $T$ in such a way that $B_T \sim \mu$ and $(B_{t \wedge T})_{t \geq 0}$ is a uniformly integrable martingale. We will see an explicit and more general construction using this approach by Bertoin and Le Jan in Section 3.15. To a lesser extent also works of Jeulin and Yor [60] and Vallois [118] can be seen as continuations of Monroe's ideas.

### 3.7. Heath (1974)

Heath [53] considers the Skorokhod embedding problem for the $n$-dimensional Brownian motion. More precisely, he considers Brownian motion killed at the first exit time from the unit ball. The construction is an adaptation of a potential theoretic construction of Mokobodzki, and relies on results of Rost (see Section 3.5). The results generalize to processes satisfying the "right hypotheses". We note that the stopping times are randomized.

### 3.8. Chacon and Walsh (1974-1976)

Making use of potential theory on the real line, Chacon and Walsh [19] gave an elegant and simple description of a solution which happens to be quite general. Their work was based on an earlier paper of Chacon and Baxter [5], who worked with a more general setup and obtained results, for example, for $n$-dimensional Brownian motion. This approach proves very fruitful in one dimension and we will try to make the Chacon and Walsh solution a reference point throughout this survey.

Recall from Definition 2.2 that, for a centered probability measure $\mu$ on $\mathbb{R}$, its potential is defined via $U\mu(x) = -\int_{\mathbb{R}} |x - y| d\mu(y)$. We will now use the properties of this functional, given in Proposition 2.3, to describe the solution proposed by Chacon and Walsh.

Write $\mu_0 = \delta_0$. Choose a point between the graphs of $U\mu_0$ and $U\mu$, and draw a line $l_1$ through this point which stays above $U\mu$ (actually in the original construction tangent lines are considered, which is natural but not necessary). This line cuts the potential $U\delta_0$ in two points $a_1 < 0 < b_1$. We consider the new potential $U\mu_1$ given by $U\mu_0$ on $(-\infty, a_1] \cup [b_1, \infty)$ and linear on $[a_1, b_1]$. We iterate the procedure. The choice of lines which we use to produce potentials $U\mu_n$ is not important. It suffices to see that we can indeed choose lines in such a way that $U\mu_n \to U\mu$ (and therefore $\mu_n \Rightarrow \mu$). This is true, as $U\mu$ is a concave function and it can be represented as the infimum of a countable number of affine functions. The stopping time is obtained therefore through a limit procedure. If we write $T_{a,b}$ for the exit time of $[a,b]$ and $\theta$ for the standard shift operator, then $T_1 = T_{a_1,b_1}$, $T_2 = T_1 + T_{a_2,b_2} \circ \theta_{T_1}$, ..., $T_n = T_{n-1} + T_{a_n,b_n} \circ \theta_{T_{n-1}}$ and $T = \lim T_n$. It is fairly easy to show that the limit is almost surely finite and for a measure with finite second moment, $\mathbb{E}T = \int x^2 d\mu(x)$ (via *(v)* in Proposition 2.3). The solution is easily explained with a series of drawings:

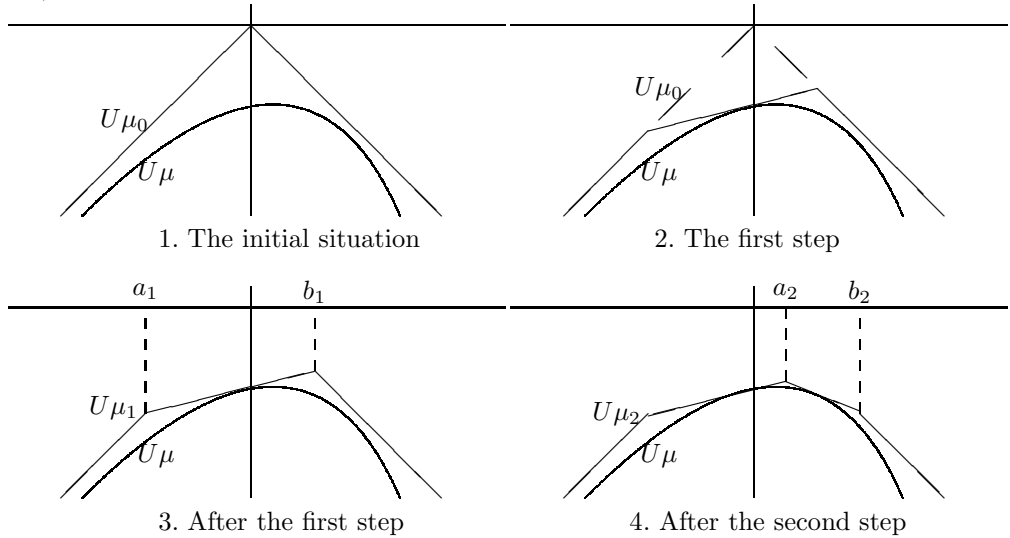

1. The initial situation
2. The first step
3. After the first step
4. After the second step

What is really being done on the level of Brownian motion? The procedure was summarized by Meilijson [73]. Given a centered measure $\mu$ express it as a limit of $\mu_n$ - measures with finite supports. Do it in such a way that there exists a martingale $(X_n)_{n\geq 1}$, $X_n \sim \mu_n$, which has almost surely dichotomous transitions (i.e. conditionally on the past with $X_n = a$, there are only two possible values for $X_{n+1}$). Then embed the martingale $X_n$ in Brownian motion by successive first hitting times of one of the two points.

Dubins' solution (see Section 3.2) is a special case of this procedure. What Dubins proposes is actually a simple method of choosing the tangents. To obtain the potential $U\mu_1$ draw tangent at $(0, U\mu(0))$, which will cut the potential $U\delta_0$ in two points $a < 0 < b$. Then draw tangents in $(a, U\mu(a))$ and $(b, U\mu(b))$. The lines will cut the potential $U\mu_1$ in four points yielding the potential $U\mu_2$. Draw the tangents in those four points obtaining 8 intersections with $U\mu_2$ which give new coordinates for drawing tangents. Then, iterate.

Chacon and Walsh solution remains true in a more general setup. It is very important to note that the only property of Brownian motion we have used here is that for any $a < 0 < b$, $\mathbb{P}(B_{T_{a,b}} = a) = \frac{b}{b-a}$. This is the property which makes the potentials of $\mu_n$ piece-wise linear, or more generally which makes the assertion *(vii)* of Proposition 2.3 true. However, this is true not only for Brownian motion but for any continuous martingale $(M_t)$ with $\langle M, M \rangle_\infty = \infty$ a.s. The solution presented here is therefore valid in this more general situation.

The methodology presented above, allows us to recover intuitively other solutions and we will make of it a reference point for us. In the case of a measure $\mu$ with finite support, this solution contains as a special case not only the solution of Dubins but also the solution of Azéma and Yor (see Section 4). The general Azéma-Yor solution can also be explained via this methodology (cf. Section 5.1) as well as the Vallois construction (cf. Section 3.12) as proven by Cox [22].

### 3.9. Azéma-Yor (1979)

The solution developed by Azéma and Yor [3] has received a lot of attention in the literature and its properties were investigated in detail. In consequence we might not be able to present all of the relative material. We will first discuss the solution for measures with finite support in Section 4, and then the general case in Section 5. The general stopping time is given in (5.3) and the Hardy-Littlewood function, on which it relies, is displayed in (5.2).

### 3.10. Falkner (1980)

Falkner [36] discusses embeddings of measures into the $n$-dimensional Brownian motion. This was then continued for right processes by Falkner and Fitzsimmons [37]. Even though we do not intend to discuss this topic here we point out some difficulties which arise in the multi-dimensional case. They are mainly due to the fact that the Brownian motion does not visit points any more. If we consider measures concentrated on the unit circle it is quite easy to see that only one of them, namely the uniform distribution, can be embedded by means of an integrable stopping time. Distributions with atoms cannot be embedded at all.

### 3.11. Bass (1983)

The solution proposed by Bass is basically an application of Itô's formula and the Dambis-Dubins-Schwarz theorem. Bass uses also some older results of Yershov [123]. This solution is also reported in Stroock's book ([117], p. 213–217). Let $F$ be the distribution function of a centered measure $\mu$, $\Phi$ the distribution function of $\mathcal{N}(0,1)$ and $p_t(x)$ the density of $\mathcal{N}(0,t)$. Define function $g(x) = F^{-1}(\Phi(x))$. Take a Brownian motion $(B_t)$, then $g(B_1) \sim \mu$. Using Itô's formula (or Clark's formula) we can write

$$g(B_1) = \int_0^1 a(s, B_s)dB_s, \text{ where } a(s,y) = -\int \frac{\partial p_{1-s}(z)}{\partial z} g(z+y)dz.$$

We put $a(s,y) = 1$ for $s \geq 1$ and define a local martingale $M$ through $M_t = \int_0^t a(s, B_s)dB_s$. Its quadratic variation is given by $R(t) = \int_0^t a^2(s, B_s)ds$. Denote by $S$ the inverse of $R$. Then the process $N_u = M_{S(u)}$, $u \geq 0$, is a Brownian motion. Furthermore $N_{R(1)} = M_1 \sim \mu$. Bass [4] then uses a differential equation argument (following works by Yershov [123]) to show that actually $R(1)$ is a stopping time in the natural filtration of $N$. The argument allows to construct, for an arbitrary Brownian motion $(\beta_t)$, a stopping time $T$ such that $(N_{R(1)}, R(1)) \sim (\beta_T, T)$, which is the desired solution to the embedding problem.

In general the martingale $M$ is hard to describe. However in simple cases it is straightforward. We give two examples. First consider $\mu = \frac{1}{2}(\delta_{-1} + \delta_1)$. Then

$$M_t = sgn(B_t)\Big(1 - 2\Phi\Big(\frac{-|B_t|}{\sqrt{1-t}}\Big)\Big), \; t \leq 1.$$

Write $g_1$ for "the last 0 before time 1," that is: $g_1 = \sup\{u < 1 : B_u = 0\}$. This is not a stopping time. The associated Azéma supermartingale (cf. Yor [124] p. 41-43) is denoted $Z_t^{g_1} = \mathbb{P}(g_1 > t|\mathcal{F}_t)$, where $(\mathcal{F}_t)$ is the natural filtration of $(B_t)$. We can rewrite $M$ in the following way

$$M_t = sgn(B_t)(1 - 2Z_t^{g_1}) = \mathbb{E}\Big[sgn(B_t)(1 - 2\mathbf{1}_{g_1>t})|\mathcal{F}_t\Big].$$

Consider now $\mu = \frac{3}{4}\delta_0 + \frac{1}{8}(\delta_{-2} + \delta_2)$ and write $\xi = -\Phi^{-1}(\frac{1}{8})$. Then

$$M_t = sgn(B_t)2\Big(\Phi\Big(\frac{-\xi + |B_t|}{\sqrt{1-t}}\Big) - \Phi\Big(\frac{-\xi - |B_t|}{\sqrt{1-t}}\Big)\Big) = \mathbb{E}\Big[\mathbf{1}_{B_1>\xi} - \mathbf{1}_{B_1<-\xi}|\mathcal{F}_t\Big].$$

What really happens in the construction of the martingale $M$ should now be clearer: we transform the standard normal distribution of $B_1$ into $\mu$ by a map $\mathbb{R} \to \mathbb{R}$ and hit this distribution by readjusting Brownian paths.

### 3.12. Vallois (1983)

In their work, Jeulin and Yor [60] studied stopping times $T^{h,k} = \inf\{t \geq 0 : B_t^+ h(L_t) + B_t^- k(L_t) = 1\}$, where $B$ and $L$ are respectively Brownian motion and its local time at 0, and $k, h$ are two positive Borel functions which satisfy some additional conditions. They described the distributions of $X_T$ and $(X_T, T)$. This, using Lévy's equivalence theorem which asserts that the two-dimensional processes $(S_t, S_t - B_t)$ and $(L_t, |B_t|)$ have the same law, can be seen as a generalization of results of Azéma and Yor [3]. Vallois [118] follows the framework of Jeulin and Yor [60] and develops a complete solution to the embedding problem for Brownian motion. For any probability measure on $\mathbb{R}$ he shows there exist two functions $h^+$ and $h^-$ and $\delta > 0$ such that

$$T_{\mathbf{V}} = T_\delta^L \wedge \inf\{t : B_t^+ = h^+(L_t) \text{ or } B_t^- = h_-(L_t)\},$$

embeds $\mu$ in $B$, $B_{T_\mathbf{V}} \sim \mu$, where $T_\delta^L = \inf\{t : L_t = \delta\}$. His method allows us to calculate the functions $h^+$ and $h^-$ as is shown in several examples. Furthermore, Vallois proves that $(B_{t \wedge T_\mathbf{V}} : t \geq 0)$ is a uniformly integrable martingale if and only if $\mu$ has a first moment and is centered.

The formulae obtained by Vallois are quite complicated. In the case of symmetric measures they can simplified considerably as noted in our paper [84].

### 3.13. Perkins (1985)

Perkins [90] investigated the problem of $H^1$-embedding and developed a solution to the Skorokhod problem that minimizes stochastically the distribution of the maximum of the process (up to the stopping time) and maximizes stochastically the distribution of the minimum. We will present his work in Section 6.1.

### 3.14. Jacka (1988)

Jacka [58] was interested in what can be seen as a converse to Perkins' goal. He looked for a maximal $H^1$-embedding, that is an embedding which maximizes stochastically the distribution of the maximum of the absolute value of the process. We describe his solution in Section 6.2. It seems that Jacka's work was mainly motivated by applications in the optimal stopping theory. In his subsequent paper, Jacka [59] used his embedding to recover the best constants $C_p$ appearing in

$$\mathbb{E}\left[\sup_{t \geq 0} M_t\right] \leq C_p \Big(\mathbb{E}\, |M_\infty|^p\Big)^{\frac{1}{p}} \quad (p > 1). \tag{3.5}$$

The interplay between Skorokhod embeddings and optimal stopping theory, of which we have an example above, is quite rich and diverse. We will come back to this matter in Section 10.

### 3.15. Bertoin and Le Jan (1993)

Tools used by Bertoin and Le Jan in [6] to construct their solution to the Skorokhod problem are similar to those employed by Monroe [79]. The authors however do not rely on Monroe's work. They consider a much more general setting and obtain explicit formulae.

Let $E$ be a locally compact space with a countable basis, and $X$ a Hunt process on $E$, with a regular recurrent starting point $0 \in E$. The local time of $X$ at $0$ is denoted by $L$. Authors use the correspondence between positive, continuous, additive functionals and their Revuz measures (see Revuz [98]). The invariant measure $\lambda$ which is used to define the correspondence is given by

$$\int f d\lambda = \int d\eta \int_0^\zeta f(X_s) ds + cf(0),$$

where $\eta$ is the characteristic measure of the excursion process of $X$, $c$ is the delay coefficient[4] and $\zeta$ is the first hitting time of 0.

Take a Revuz measure $\mu$ with $\mu(\{0\}) = 0$, $\mu(E) \leq 1$ and $A^\mu$ the positive continuous additive functional associated with it. Then $T = \inf\{t : A_t^\mu > L_t\}$ embeds $\mu$ in $X$. However this solution is not a very useful one as we have $\mathbb{E} L_T = \infty$.

Consider an additional hypothesis:

*There exists a bounded Borel function $\hat{V}_\mu$, such that for every Revuz measure $\nu$:* $\int \hat{V}_\mu d\nu = \int \mathbb{E}_x(A_\zeta^\nu) d\mu(x).$ (3.6)

Then, if $\mu$ is a probability measure and $\gamma_0 = ||\hat{V}_\mu||_\infty$, for any $\gamma \geq \gamma_0$, the following stopping time embeds $\mu$ in $X$:

$$T_{\mathbf{BLJ}}^\gamma = \inf\left\{t \geq 0 : \gamma \int_0^t (\gamma - \hat{V}_\mu(X_s))^{-1} dA_s^\mu > L_t\right\}.$$ (3.7)

Furthermore this solution is optimal in the following sense. Consider any Revuz measure $\nu$ with $A^\nu$ not identically 0. Then $\mathbb{E}_0 A_{T_{\mathbf{BLJ}}^\gamma}^\nu = \int (\gamma - \hat{V}_\mu) d\nu$ and for any solution $S$ of the Skorokhod problem for $\mu$, $\mathbb{E}_0 A_S^\nu \geq \int (\gamma_0 - \hat{V}_\mu) d\nu$.

Authors then specify to the case of symmetric Lévy processes. Let us investigate this solution in the case of $X$ a real Brownian motion. The function $\hat{V}_\mu$ exists and its supremum is given by $\gamma_0 = ||\hat{V}_\mu||_\infty = 2\max\{\int x^+ d\mu, \int x^- d\mu\}$. Hypothesis 3.6 is equivalent to $\mu$ having a finite first moment. In this case introduce

$$\rho(x) = \begin{cases} 2\int_{a>x}(a-x)d\mu(a) & , \text{ for } x \geq 0 \\ 2\int_{a<x}(x-a)d\mu(a) & , \text{ for } x < 0. \end{cases}$$

Then

$$T_{\mathbf{BLJ}}^{\gamma_0} = \inf\left\{t : \gamma_0 \int \frac{L_t^x}{\rho(x)} d\mu(x) > L_t^0\right\}.$$

Note that if $supp(\mu) \subset [a,b]$ then $\rho(a) = \rho(b) = 0$, so in particular $T_{\mathbf{BLJ}} \leq T_{a,b}$, the first exit time from $[a,b]$.

### 3.16. Vallois (1994)

In Section 3.12 above we described a solution to the Skorokhod embedding problem, based on Brownian local time, developed by Vallois in [118]. Here we present a different solution which Vallois described in [120]. Actually, Vallois considered a generalization of the Skorokhod embedding problem.

Let $(M_t : t \geq 0)$ be a continuous local martingale with $M_\infty \sim \mu$. Vallois [118] gave a complete characterization of the class of possible laws of $S_\infty = \sup_{t\geq 0} M_t$ (see also Section 10.3). To prove his characterization he had to solve the following embedding problem: let $\mu \in \mathcal{M}_0^1$ and $\mu_1$ be a sub-probability measure on $\mathbb{R}_+$ such that $\mu - \mu_1$ is a positive measure on $\mathbb{R}_+$. For $(B_t : t \geq 0)$ a Brownian motion, and $S_t = \sup_{u \leq t} B_u$, construct a stopping time $T$ such that $B_T \sim \mu$, $B_T|_{\{B_T = S_T\}} \sim \mu_1$ and $(B_{t \wedge T} : t \geq 0)$ is a uniformly integrable martingale. Vallois [120], showed that for each such couple $(\mu, \mu_1)$ there exists an increasing function $\alpha$ such that the stopping time

$$T_{\mu,\mu_1} = \inf\{t \geq 0 : (B_t = \alpha(S_t), S_t \notin \Gamma) \text{ or } B_t = A\}, \qquad (3.8)$$

where $\Gamma = \{x \in supp(\mu) : \alpha(x) = x\}$ and $A$ is an independent random variable whose distribution is a specified function of $\mu$ and $\mu_1$, solves the embedding problem under consideration.

---

[4]That is: the holding time at 0 has an exponential distribution with parameter $c$.

Vallois' construction provides a link between two extreme cases: $\mu_1 = 0$ and $\mu_1 = \mu|_{[0,\infty)}$. In the first case, $\mu_1 = 0$, we just have $T_{\mu,0} = T_{\mathbf{AY}}$, where $T_{\mathbf{AY}}$ is the Azéma-Yor stopping time defined via (5.3). This case yields the maximal possible distribution of $S_T$ (see Section 5.5). The second case, $\mu_1 = \mu|_{[0,\infty)}$, yields the minimal possible distribution of $S_T$, which coincides with $S_{T_{\mathbf{P}}}$, where Perkins' solution $T_{\mathbf{P}}$ is given by (6.3).

We believe that one could see Vallois' reasoning via one-dimensional potential theory but we were not able to do so. This is partially due to the fact that the stopping rule obtained is randomized. Vallois [120] in his proof used martingale theory, in particular he relied on the martingales displayed in (5.11).

### 3.17. Falkner, Grandits, Pedersen and Peskir (2000 - 2002)

We mention, exceptionally, several solutions in this one subsection, the reason being that they are all closely connected and subsequent ones were thought of as generalizations of the preceding. The solution developed by Azéma and Yor [3] works for recurrent real-valued diffusions. It is natural to ask whether this solution could be generalized for transient diffusions.

Brownian motion with drift was the first transient diffusion considered. It was done by Hall back in 1969 as we stressed in Section 3.4, however his work doesn't seem to be well known. Consequently, the subject was treated again in a similar manner by Grandits and Falkner [46] and Peskir [93]. They show that for a probability measure $\mu$ there exists a stopping time $T$ such that $B_T^{(\delta)} = B_T + \delta T \sim \mu$ ($\delta > 0$) if and only if $\int e^{-2\delta x} d\mu(x) \leq 1$. The authors work with the martingale $Y_t = \exp(-2\delta B_t^{(\delta)}) - 1$ rather than with the process $B^{(\delta)}$ itself and they adapt the arguments used by Azéma and Yor. Grandits and Falkner gave an explicit formula for $T$ when $\int e^{-2\delta x} d\mu(x) = 1$, while Peskir obtained an explicit formula valid in all cases and showed that the solution maximizes stochastically the supremum of the process up to the stopping time. He also showed that letting $\delta \to 0$ one recovers the standard Azéma-Yor solution. The arguments used were extended by Pedersen and Peskir [89], to treat any nonsingular diffusion on $\mathbb{R}$ starting at 0.

The methodology of these embeddings can be summarized as follows. Given a diffusion, compose it with its scale function. The result is a continuous local martingale. Use the Dambis-Dubins-Schwarz theorem to work with a Brownian motion. Then use the Azéma-Yor embedding and follow the above steps in a reversed direction. Some care is just needed as the local martingale obtained from the diffusion might converge. We develop this subject in Section 9.

### 3.18. Roynette, Vallois and Yor (2002)

These authors [108] consider Brownian motion together with its local time at 0 as a Markov process $X$, and apply Rost's results (see Section 3.5) to derive conditions on $\mu$ - a probability measure on $\mathbb{R} \times \mathbb{R}_+$ - under which there exists an embedding of $\mu$ in $X$. They give an explicit form of the condition and of the solution whenever there exists one.

Consider $\mu(dx, dl) = \nu(dl)M(l, dx)$ a probability measure on $\mathbb{R} \times \mathbb{R}_+$, where $\nu$ is a probability measure on $\mathbb{R}_+$ and $M$ is a Markov kernel. Then there exists a stopping time $T_{\mathbf{RVY}}$, such that $(B_{T_{\mathbf{RVY}}}, L_{T_{\mathbf{RVY}}}) \sim \mu$ and $(B_{t \wedge T_{\mathbf{RVY}}} : t \geq 0)$ is a uniformly integrable martingale if and only if

1. $\int_{\mathbb{R} \times \mathbb{R}_+} |x| \mu(dx, dl) < \infty$,

2. $m(l) := \int_{\mathbb{R}} x^+ M(l, dx) = \int_{\mathbb{R}} x^- M(l, dx) < \infty$,

3. $m(l)\nu(dl) = \frac{1}{2}\nu[l, \infty) dl$.

Moreover, if $\nu$ is absolutely continuous with respect to Lebesgue measure, then $T$ is a stopping time with respect to $\mathcal{F}_t = \sigma(B_s : s \leq t)$.

An explicit construction of $T_{\mathbf{RVY}}$ is given. It is carried out in two steps. First, a stopping time $T_1$ is constructed which embeds $\nu$ in $L$ - it follows arguments of Jeulin and Yor [60]. Then, using the Azéma-Yor embedding, $T > T_1$ is constructed such that conditionally on $L_{T_1} = l$, $B_T \sim M(l, dx)$ and $B_t$ has no zeros between $T_1$ and $T$, so that $L_{T_1} = L_T$.

This construction yields two interesting studies. It enables the authors in [108] to investigate the possible laws of $L_T$, where $T$ is a stopping time such that $(B_{t \wedge T} : t \geq 0)$ is a uniformly integrable martingale, thus generalizing the work of Vallois [120]. It also gives the possibility to look at the question of independence of $T$ and $B_T$ (primarily undertaken in [28], [29]).

### 3.19. Hambly, Kersting and Kyprianou (2002)

These authors [51] develop a Law of Iterated Logarithm for random walks conditioned to stay nonnegative. One of the ingredients of their proof is an embedding of a given distribution into a Bessel(3) process. This distribution is actually the distribution of a renewal function of the initial distribution under the law of random walk conditioned to stay nonnegative. Authors use properties of the latter to realize their embedding. It is done in two steps, as they first stop the process in an intermediate distribution. Both parts of the construction require external randomization, independent of the process.

### 3.20. Cox and Hobson (2004)

The work of Cox and Hobson [23] can be seen as a generalization of Perkins' solution (see Section 6.1) to the case of regular diffusions and measures which are not necessarily centered. Indeed, consider any probability measure $\mu$ on $\mathbb{R}$. The authors define two functions $\gamma_-$ and $\gamma_+$, which for a centered probability measure $\mu$ coincide with the ones defined by Perkins. For a continuous local martingale $(M_t)$ define its maximum and minimum processes through $S_t^M = \sup_{s \leq t} M_s$ and $s_t^M = -\inf_{s \leq t} M_s$. Then, given that $\langle M, M \rangle_\infty = \infty$, the stopping time $T_{\mathbf{CH}} = \inf\{t > 0 : M_t \notin (-\gamma_-(S_t^M), \gamma_+(s_t^M))\}$ embeds $\mu$ in $M$. Furthermore, it minimizes stochastically the supremum and maximizes stochastically the minimum (cf. (6.4), (6.5)). The setup is then extended to regular diffusions using scale functions. For regular diffusions, scale functions exist and are continuous, strictly increasing. This implies that the maximum and minimum are preserved. Conditions on the diffusion (through its scale function) and on $\mu$, for the existence of an embedding, are given.

Further, Cox and Hobson investigate the existence of $H^p$-embeddings (see Section 6). Conditions are given which in important cases are necessary and sufficient.

### 3.21. Obłój and Yor (2004)

Almost all explicit solutions discussed so far dealt with continuous processes. The solution proposed by Bertoin and Le Jan (cf. Section 3.15) is the only exception, and even this solution is explicit only if we have easy access to additive functionals and can identify the function $\hat{V}_\mu$ given by (3.6). Stopping discontinuous processes is in general harder, as we have less control over the value at the stopping time. Yet, in some cases, we can develop tools in discontinuous setups inspired by the classical framework. This is the case of an embedding for "nice" functionals of Brownian excursions, described in our joint work with Yor [84]. This solution was inspired by the Azéma-Yor solution (see Section 5), as we explain below. For the sake of simplicity we will not present the reasoning found in Obłój and Yor [84] in all generality but we will rather restrain ourselves to the case of the age process. We are interested therefore in developing an embedding for the following process: $A_t = t - g_t$, where $g_t = \sup\{s \leq t : B_s = 0\}$ is the last zero before time $t$. In other words, $A_t$ is the age, at time $t$, of the Brownian excursion straddling time $t$. This is a process with values in $\mathbb{R}_+$, which is discontinuous and jumps down to 0 at points of discontinuity. More precisely, we will focus on an equivalent process given by $\tilde{A}_t = \sqrt{\frac{\pi}{2}(t - g_t)}$. This process is, in some sense, more natural to consider, as it is just the projection of the absolute value of Brownian motion on the filtration generated by the signs.

Starting with the standard Azéma-Yor embedding (see Section 5.1), let us rewrite their stopping time using Lévy's equivalence theorem.

$$T_{\mathbf{AY}} \stackrel{\text{see (5.3)}}{=} \inf\{t \geq 0 : S_t \leq \Psi_\mu(B_t)\} = \inf\{t \geq 0 : S_t - \Psi_\mu^{-1}(S_t) \geq S_t - B_t\}$$

$$\stackrel{law}{=} \inf\{t \geq 0 : \tilde{\Psi}_\mu(L_t) \geq |B_t|\}, \quad \text{where} \quad \tilde{\Psi}_\mu = Id - \Psi_\mu^{-1}. \tag{3.9}$$

The idea of replacing the supremum process by the local time process comes from the fact that the latter, unlike the former, is measurable with respect to the filtration generated by the signs (or by the zeros) of Brownian motion. Even though the function $\tilde{\Psi}_\mu$ is in general hard to describe, the above formula suggests that it might be possible to develop an embedding for the age process of the following form: $T = \inf\{t > 0 : \tilde{A}_t \geq \varphi(L_t)\}$, where the function $\varphi$ would be a function of the measure $\mu$ we want to embed. This proves to be true. For any measure $\mu$ on $\mathbb{R}_+$ with $\mu(\{0\}) = 0$, let $[a,b]$ be its support, $\overline{\mu}(x) = \mu([x,\infty))$ its tail, and define the dual Hardy-Littlewood function[5]

$$\psi_\mu(x) = \int_{[0,x]} \frac{y}{\overline{\mu}(y)} d\mu(y), \quad \text{for } a \leq x < b, \tag{3.10}$$

$\psi_\mu(x) = 0$ for $0 \leq x < a$ and $\psi_\mu(x) = \infty$ for $x \geq b$. The function $\psi_\mu$ is right-continuous, non-decreasing and $\psi_\mu(x) < \infty$ for all $x < b$. We can then define its right inverse $\varphi_\mu = \psi_\mu^{-1}$. Then the stopping time

$$T_{\mathbf{OY}} = \inf\{t > 0 : \tilde{A}_t \geq \varphi_\mu(L_t)\}$$

embeds $\mu$ in $\tilde{A}$, $\tilde{A}_{T_{\mathbf{OY}}} \sim \mu$.

The proof is carried out using excursions theory (even though a different one, relying on martingale arguments, is possible too). The crucial step is the calculation of the law of $L_{T_{\mathbf{OY}}}$, which we give here. Note $(\tau_l)$ the inverse of the local time at 0, $(e_l)$ the excursion process and $V(e_l) = \tilde{A}_{\tau_l-}$ the rescaled lifetime of the excursion.

$$\begin{aligned}
\mathbb{P}(L_{T_{\mathbf{OY}}} \geq l) &= \mathbb{P}(T_{\mathbf{OY}} \geq \tau_l) \\
&= \mathbb{P}\Big(\text{on the time interval } [0, \tau_l] \text{ for every excursion } e_s,\ s \leq l, \\
&\qquad \text{its (rescaled) lifetime did not exceed } \varphi_\mu(s)\Big) \\
&= \mathbb{P}\Big(\sum_{s \leq l} \mathbf{1}_{V(e_s) \geq \varphi_\mu(s)} = 0\Big) = \exp\Big(-\int_0^l \frac{ds}{\varphi_\mu(s)}\Big),
\end{aligned} \tag{3.11}$$

where the last equality follows from the fact that the process $N_l = \sum_{s \leq l} \mathbf{1}_{F(e_s) \geq \varphi(s)}$, $l \geq 0$, is an inhomogeneous Poisson process with appropriate parameter. The core of the above argument remains valid for much more general functionals of Brownian excursions, thus establishing a methodology of solving the Skorokhod embedding in a number of important setups. In particular, it yields an embedding for the absolute value $(|B_t| : t \geq 0)$, giving the explicit formula in the Vallois solution (see Section 3.12) for symmetric measures.

---

[5]Compare with the Hardy-Littlewood function $\Psi_\mu$ used by Azéma and Yor given in (5.2).

| Solution | Further developments |
|---|---|
| Skorokhod (1961) [115] | Strassen [116], Breiman [11], Perkins [90] Meilijson [74] |
| Dubins (1968) [31] | Meyer [76], Bretagnolle [12] Chacon and Walsh [19] |
| Hall, W.J. (1968) [49] | |
| Root (1969) [103] | Loynes [70], Rost [107] |
| Rost (1971) [105] | Azéma and Meyer [1], Baxter, Chacon P. [16], Fitzsimmons [38], [39], [40] |
| Monroe (1972) [79] | Vallois [118], Bertoin and Le Jan [6] |
| Heath (1974) [105] | |
| Chacon and Walsh (1974-76) [5], [19] | Chacon and Ghoussoub [18] |
| Azéma and Yor (1979) [3] | Azéma and Yor [2], Pierre [95], Jeulin and Yor [60] Rogers [101], Meilijson [73], Zaremba [126] van der Vecht [121], Pedersen and Peskir [88] |
| Falkner (1980) [36] | Falkner and Fitzsimmons [37] |
| Bass (1983) [4] | |
| Vallois (1983) [118] | |
| Perkins (1985) [90] | Cox and Hobson [23] |
| Jacka (1988) [58] | |
| Bertoin and Le Jan (1992) [6] | |
| Grandits and Falkner [46], Peskir (2000) [93] | Pedersen and Peskir [89] |
| Pedersen and Peskir (2001) [89] | Cox and Hobson [23] |
| Roynette, Vallois and Yor (2002) [108] | |
| Hambly, Kersting and Kyprianou (2002) [51] | |
| Cox and Hobson (2004) [23] | Cox and Hobson [24] |
| Obłój and Yor (2004) [84] | |

Table 1: Genealogy of solutions to the Skorokhod embedding problem

## 4. Embedding measures with finite support

We continue here the discussion of Chacon and Walsh's solution from Section 3.8. The target measure $\mu$, which is a centered probability measure on $\mathbb{R}$, is supposed to have a finite support, i.e.: $\mu = \sum_{i=1}^{n} \alpha_i \delta_{x_i}$, where $\alpha_i > 0$, $\sum_{i=1}^{n} \alpha_i = 1$ and $\sum_{i=1}^{n} \alpha_i x_i = 0$, $x_i \neq x_j$ for $1 \leq i \neq j \leq n$. Recall the potential $U\mu$ given in Definition 2.2. In the present situation, it is a piece-wise linear function breaking at the atoms $\{x_1, \ldots, x_n\}$. It is interesting that in this simple setting several solutions are easily seen to be special cases of the Chacon-Walsh solution.

### 4.1. Continuous local martingale setup

Let $(M_t)$ be a continuous local martingale, $M_0 = 0$, with $\langle M, M \rangle_\infty = \infty$ a.s. and consider the embedding problem of $\mu$ in $M$. First we can ask: how many solutions to the embedding problem exist? The answer is simple: except for the trivial case, an infinite number. Recall that $n$ is the number of points in the support of $\mu$. We say that two solutions $T$ and $S$ that embed $\mu$ are different if $\mathbb{P}(T = S) < 1$.

**Proposition 4.1.** *If $n \geq 3$, then there exists uncountably many stopping times $T$ which embed $\mu$: $M_T \sim \mu$, and such that $(M_{t \wedge T})$ is a uniformly integrable martingale.*

*Proof.* This is an easy consequence of the Chacon and Walsh embedding described in Section 3.8. We describe only how to obtain an infinity of solutions for $n = 3$, the ideas for higher values of $n$ being the same. Our reasoning is easiest explained with a series of drawings:

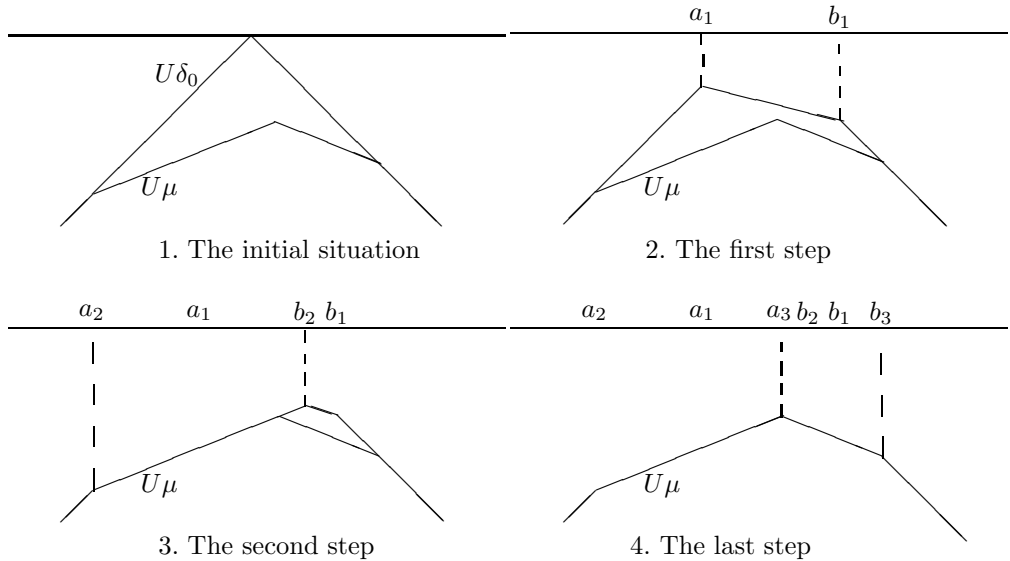

Note that $x_1 = a_2$, $x_2 = a_3$ and $x_3 = b_3$. Consider the set of paths $\omega$ which are stopped in $a_2 = x_1$: $A_T = \{\omega : M(\omega)_{T(\omega)} = x_1\}$. It can be described as the set of paths that hit $a_1$ before $b_1$ and then $a_2$ before $b_2$. It is clear from our construction that we can choose these points in an infinite number of ways and each of them yields a different embedding (paths stopped at $x_1$ differ). Indeed as for any $a_1 \leq \overline{a}_1 < 0 < b_2 < \overline{b}_2 < b_1$, $\mathbb{P}(T_{\overline{a}_1} < T_{\overline{b}_2} < T_{a_1} < T_{b_1}) > 0$, the result follows from the fact that if $\mathbb{P}(B_T \neq B_S) > 0$ then $\mathbb{P}(S = T) < 1$. □

However in their original work Chacon and Walsh allowed only to take lines tangent to the graph of $U\mu$. In this case we have to choose the order of drawing $n-1$ tangent lines and this can be done in $(n-1)!$ ways. Several of them are known as particular solutions. Note that we could just as well take a measure with a countable support and the constructions described below remain true.

- **The Azéma-Yor** solution consists in drawing the lines from left to right. This was first noticed, for measures with finite support, by Meilijson [73]. The Azéma-Yor solution can be summarized in the following way. Let $X \sim \mu$. Consider a martingale $Y_i = \mathbb{E}[X|\mathbf{1}_{X=x_1}, \ldots, \mathbf{1}_{X=x_i}]$, $Y_0 = 0$. Note that $Y_i$ has dichotomous transitions, it increases and then becomes constant and equal to $X$ as soon as it decreases, $Y_{n-1} = X$. Embed this martingale into $(M_t)$ by consecutive first hitting times. It suffices to remark now that the hitting times correspond to ones obtained by taking tangent lines in Chacon and Walsh construction, from left to right. Indeed, the lower bounds for the consecutive hitting times are just atoms of $\mu$ in increasing order from $x_1$ to $x_{n-1}$.

  Another way to see all this is to use the property that the Azéma-Yor solution maximizes stochastically the maximum of the process. If we want the supremum of the process to be maximal we have to take care not to "penalize big values" or in other words "not to stop too early in a small atom" - this is just to say that we want to take the bounds of our consecutive stopping times as small as possible, which is just to say we take the tangent lines from left to right.

  Section 5 is devoted to the study of this particular solution in a general setting.

- **The Reversed Azéma-Yor** solution consists in taking tangents from right to left. It minimizes the infimum of the stopped process. We discuss it in Section 5.3.

- **The Dubins** solution is always a particular case of the Chacon and Walsh solution. In a sense it gives a canonical method of choosing the tangents. We described it in Section 3.8.

### *4.2. Random walk setup*

When presenting above, in Section 3.8, the Chacon and Walsh solution, we pointed out that it only takes advantage of the martingale property of Brownian motion. We could ask therefore if it can also be used to construct embeddings for processes with countable state space. And it is very easy to see the answer is affirmative, under some hypothesis which guarantees that no over-shoot phenomenon can take place.

Let $(M_n)_{n\geq 0}$ be a discrete martingale taking values in $\mathbb{Z}$ with $M_0 = 0$ a.s. Naturally, a symmetric random walk provides us with a canonical example and actually this section was motivated by a question raised by Fujita [44] concerning precisely the extension of the Azéma-Yor embedding to random walks.

Denote the range of $M$ by $R = \{z \in \mathbb{Z} : \inf\{n \geq 0 : M_n = z\} < \infty\ a.s.\}$ and let $\mu$ be a centered probability measure on $R$. It is then an atomic measure[6]. We can still use the Chacon and Walsh methodology under one condition. Namely, we can only use such exit times $T_{a,b}$ that both $a, b \in R$. Otherwise we would face an over-shoot phenomena and the method would not work.

Naturally in the Azéma-Yor method, the left end points, $(a_k)$, are just the consecutive atoms of $\mu$, which by definition are in $R$. The right-end points are determined by tangents to $U\mu$ as is described in detail in Section 5.1 below. They are given by $\Psi_\mu(a_k)$, where $\Psi_\mu$ is displayed in (5.2), and they might not be in $R$.

**Proposition 4.2.** *In the above setting suppose that the measure $\mu$ is such that $\Psi_\mu(\mathbb{R}) \subset R$, where $\Psi_\mu$ is given by (5.2). Then the Azéma-Yor stopping time, $T = \inf\{n \geq 0 : S_n \geq \Psi_\mu(M_n)\}$, where $S_n = \max_{k \leq n} M_k$, embeds the measure $\mu$ in $M$, $M_T \sim \mu$.*

The same Proposition, for a random walk, was also obtained by Fujita [44]. He used discrete martingales, which are closely related to the martingales used by Azéma and Yor [3] in their original proof, and which we discuss in Section 5.4.

We close this section with two observations. First, notice that the condition $\Psi_\mu(\mathbb{R}) \subset R$ is not necessary for the existence of a solution to the Skorokhod embedding problem. We give a simple example. Consider a standard, symmetric random walk, $(X_n : n = 0, 1, \ldots)$, and take $\mu = \frac{2}{9}\delta_{-3} + \frac{4}{9}\delta_0 + \frac{1}{3}\delta_2$. Then $\Psi_\mu(0) = \frac{6}{7} \notin \mathbb{Z}$, but one can easily point out a good embedding. It suffices to take tangents to $U\mu$ from right to left, that is to consider the reversed Azéma-Yor solution. Then the stopping time $T = \inf\{n \geq T_{-1,2} : M_n \notin [-3, 0]\}$ embeds $\mu$ in $M$.

Finally, it is not true that for any centered discrete probability measure $\mu$, there exists a stopping time $T$, in the natural filtration of $X$, such that $X_T \sim \mu$ and $(X_{n\wedge T} : n \geq 0)$ is a uniformly integrable martingale. This was observed by Cox [21]. It suffices to take $\mu = \frac{1}{3}(\delta_{-1} + \delta_0 + \delta_1)$ (one can also impose $\mu(\{0\}) = 0$ and still find easy examples). Therefore, for some measures, a solution to the Skorokhod embedding problem for random walk has to rely on an external randomization. This subject requires further investigation.

## 5. Around the Azéma-Yor solution

In this section we present various results concerning the solution proposed by Azéma and Yor [3]. We first provide a very simple proof of the solution in

---

[6]The support is not necessary finite but this is not crucial for the method, so we decided to place this discussion in this section.

Section 5.1. The proof is based on the methodology of Chacon and Walsh (see Section 3.8) and allows us very easily to generalize the Azéma-Yor embedding for the case of nontrivial initial laws. This is done in Section 5.2. The proof encouraged us too to give, in Section 5.3, another, reversed, form of the solution. In subsequent Section 5.4, we explain the original ideas, based on martingale theory, which stood behind the work of Azéma and Yor. This allows us, in Section 5.5 to prove easily certain extremal properties of this solution.

We give first a brief history of this solution. It was presented in the $13^{th}$ "Séminaire de Probabilités" in [3] together with a follow up [2], which examined some maximizing property (see Section 5.5 below) of the solution. The original arguments are also found in Revuz and Yor ([99] pp. 269–276). The construction in [3] was carried out for a continuous local martingale and extended easily to recurrent, real-valued diffusions. The original proof was then simplified by Pierre [95] and Vallois [118]. A proof using excursion theory was given by Rogers [101]. An elementary proof was given by Zaremba [126]. This work [126] not only gives a "bare-handed" proof, but doing so it actually provides a convenient limiting procedure that allows us to see that Azéma and Yor's solution is, in a sense, a special case of Chacon and Walsh's solution (see also Section 4). The Azéma-Yor solution has been generalized by several authors - we described their papers together in Section 3.17. Finally, using Lévy's equivalence $(S-B, S) \sim (|B|, L)$, works using the local time at 0 to develop a solution to the embedding problem can be seen as generalizing Azéma and Yor's construction, and we mention therefore articles by Jeulin and Yor [60] and by Vallois [118] (see Section 3.12).

### 5.1. An even simpler proof of the Azéma-Yor formula

*In this section, for simplicity of notation, we write "AY solution" for Azéma and Yor's solution.*

We saw in Section 4 how the AY construction for a measure with finite support is just a special case of the method described by Chacon and Walsh, and consists in taking as tangent lines the extensions of the linear pieces of the potential from left to right. In this section we will see how to generalize this idea.

Let $\mu$ be a centered probability measure on $\mathbb{R}$ and $x$ a point in the interior of its support, such that $\mu$ has no atom at $x$: $\mu(\{x\}) = 0$. Recall the one-dimensional potential given in Definition 2.2. We are interested in the tangent line $l_x$ to $U\mu$ at $(x, U\mu(x))$ and its intersection with the line $\{(t, -t) : t \in \mathbb{R}\}$. Denote this intersection by $(\Psi_\mu(x), -\Psi_\mu(x))$. The following drawing summarizes these quantities:

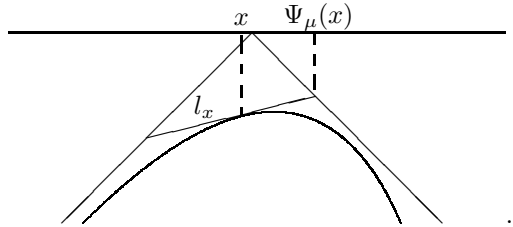

It is straightforward to find the equation satisfied by $l_x$. Its slope is given by $\frac{dU\mu(t)}{dt}|_{t=x}$ and so we can write $l_x = \{(x,y) : y = \frac{dU\mu(t)}{dt}|_{t=x} \cdot x + \beta\}$. Differentiating (2.3) we obtain

$$\begin{aligned}\frac{dU\mu(t)}{dt} &= 2\mu([t,\infty)) - 1 - 2t\mu(dt) + 2t\mu(dt) \\ &= 2\mu([t,\infty)) - 1.\end{aligned} \quad (5.1)$$

We go on to calculate $\beta$. Using the fact that $(x, U\mu(x)) \in l_x$ together with (2.3) and (5.1), we obtain $\beta = -2\int_{[x,\infty)} y d\mu(y)$. We are ready to determine the value of $\Psi_\mu(x)$ writing the equation of the intersection of $l_x$ with $\{(t,-t) : t \in \mathbb{R}\}$:

$$-\Psi_\mu(x) = \Big(2\mu([x,\infty)) - 1\Big)\Psi_\mu(x) - 2\int_{[x,\infty)} y d\mu(y), \quad \text{hence}$$

$$\Psi_\mu(x) = \frac{1}{\mu([x,\infty))} \int_{[x,\infty)} y d\mu(y), \quad (5.2)$$

which is just the barycentre function defined by Azéma and Yor, which is also called the Hardy-Littlewood (maximal) function. We still have to deal however with atoms of $\mu$. Clearly in atoms the potential is not differentiable and the tangent is not unique. It is easy to see from the construction of $\Psi_\mu(x)$ that it is a non-decreasing, continuous function (for $x$ with $\mu(\{x\}) = 0$). Fix an $x$ such that $\mu(\{x\}) > 0$. To assign value to $\Psi_\mu(x)$ we follow the case of measures with finite support described in Section 4. We see that, as a consequence of taking the tangents from left to right, we treat the point $(x, U\mu(x))$ as belonging rather to the linear piece of the potential on the left than the one on the right. In other words we take $\Psi_\mu$ to be left continuous. The formula (5.2) remains valid for all $x$ in the interior of support of $\mu$ then.

We want finally to describe the function $\Psi_\mu$ outside of the support of $\mu$ (and on its boundaries). Suppose that the support of $\mu$ is bounded from below: $a_\mu = \inf\{x : x \in supp(\mu)\} > -\infty$. For all $x \leq a_\mu$ we have simply that $U\mu(x) = x$ and so the tangent is just the line $\{(t,t) : t \in \mathbb{R}\}$, which intersects the line $\{(t,-t) : t \in \mathbb{R}\}$ in $(0,0)$, so we put $\Psi_\mu(x) = 0$.

Suppose now that the support of $\mu$ is bounded from above: $b_\mu = \sup\{x : x \in supp(\mu)\} < \infty$. Then for all $x \geq b_\mu$ we have $U\mu(x) = -x$ and so the very point $(x, U\mu(x))$ lies on the line $\{(t,-t) : t \in \mathbb{R}\}$ and we put $\Psi_\mu(x) = x$. This completes the definition of the function $\Psi_\mu$, which coincides with the one given in [3].

We want now to define the stopping time that results from this setup. At first glance we could have a small problem: following closely the case of measures with finite support would end up in a composition of an uncountable infinity of first hitting times of intervals $(x, \Psi_\mu(x))$ for $x$ in the support of $\mu$. This would not make much sense. However there is an easy way to overcome this difficulty - we can use the fact that $\Psi_\mu$ is increasing. Let $(M_t : t \geq 0)$ be a continuous local martingale and denote its supremum by $S_t = \sup\{M_s : s \leq t\}$. $S_t$ then describes which values of type $\Psi_\mu(x)$ have been reached so far, and so which of the intervals $(x, \Psi_\mu(x))$ we are considering. We want to stop when exiting such an interval on its left boundary, i.e. when we reach an $x$ such that the supremum so far is equal to $\Psi_\mu(x)$. In other words we take:

$$T_{\mathbf{AY}} = \inf\{t \geq 0 : S_t \geq \Psi_\mu(M_t)\}. \tag{5.3}$$

**Theorem 5.1 (Azéma and Yor [3]).** *Let $(M_t)_{t\geq 0}$ be a continuous martingale, $M_0 = 0$ a.s., with $\langle M, M\rangle_\infty = \infty$ a.s. For any centered probability measure $\mu$ on $\mathbb{R}$, $T_{\mathbf{AY}}$ defined through (5.3), is a stopping time in the natural filtration of $(M_t)$ and $M_{T_{\mathbf{AY}}} \sim \mu$. Furthermore, $(M_{t\wedge T_{\mathbf{AY}}})_{t\geq 0}$ is a uniformly integrable martingale.*

*Proof.* From (5.3) it is evident that $T_{\mathbf{AY}}$ is a stopping time in the natural filtration of $(M_t)$. The proof of the remaining statements is carried out through a limit procedure using approximation of $U\mu$ by piece-wise linear functions. We define measures $\mu_n$ through their potentials. Let $A = \{a_1, a_2 \ldots\}$ be the set of atoms of $\mu$ (i.e. $\mu(\{x\}) > 0$ iff $\exists! \ i : x = a_i$) and $P_n = \{\frac{k}{n} : k \in \mathbb{Z}\}$. Define $\mu_n$ through:

- $\forall x \in \{A \cup P_n\} \cap [-n, n] \quad U\mu_n(x) = U\mu(x)$,

- $U\mu_n$ is piece-wise linear, $\mu_n((-\infty, -n-1] \cup [n+1, \infty)) = 0$.

It is straightforward that $U\mu_n(x) \to U\mu(x)$ for all $x \in \mathbb{R}$, which in turn implies $\mu_n \Rightarrow \mu$. For each of the measures $\mu_n$ consider the AY stopping time $T_n$, which embeds $\mu_n$ as shown in Section 4. Finally from the definition of a tangent, we see that lines defining $T_n$ converge to tangents of $U\mu$ and $T_n \to T = T_{\mathbf{AY}}$ almost surely. Then $M_{T_n} \to M_T$ a.s. and since $\mu_n \Rightarrow \mu$, we obtain $M_{T_{\mathbf{AY}}} \sim \mu$. To see that $(M_{t\wedge T_{\mathbf{AY}}} : t \geq 0)$ one can use assertion *(v)* of Proposition 2.3 and the fact that $\mathbb{E}T_n = \int x^2 d\mu_n(x) < \infty$. $\square$

It has to be noted that what we did in this section is essentially a simplified reading of Zaremba's work [126] - simplified through the work of Chacon and Walsh.

### *5.2. Nontrivial initial laws*

One of the advantages of the approach presented here is that we impose only one condition on the initial law: namely that its potential is greater than the potential of the measure $\mu$ we want to embed (and this is known to be a necessary condition, see Section 8). We can still give an explicit formulae for the stopping time $T_{\mathbf{AY}}$, for nontrivial initial laws.

Let $(M_t : t \geq 0)$ be a continuous, local, real-valued martingale with $\langle M, M \rangle_\infty = \infty$ a.s. and with initial law $M_0 \sim \nu$, $\int |x| d\nu(x) < \infty$. Let $\mu$ be any probability measure on $\mathbb{R}$ with finite first moment, i.e. $\mu \in \mathcal{M}^1$, such that $U\mu \leq U\nu$. Using assertion *(ii)* from Proposition 2.3, this implies $\int x d\nu(x) = \int x d\mu(x) = m$, and

$$\lim_{|x| \to \infty} (U\nu(x) + |x - m|) = \lim_{|x| \to \infty} (U\mu(x) + |x - m|) = 0. \tag{5.4}$$

Then $\mu$ can be embedded in $M$ by a stopping time, which is a generalization of Azéma and Yor's solution.

Consider a tangent to the potential $U\mu$ at a point $(x, U\mu(x))$ (which is not an atom of $\mu$). We want to describe the points of intersection of this line with the potential of $\nu$. Studying (2.3), (5.1) and (5.2) it is easily seen that the points are of form $(t, U\nu(t))$, where $t$ satisfies

$$U\nu(t) = \Big(\mu([x,\infty)) - \mu((-\infty,x))\Big)t + \int_{-\infty}^{x} y d\mu(y) - \int_{x}^{\infty} y d\mu(y). \tag{5.5}$$

We want to define a generalized barycentre function $\Psi_\mu^\nu(x)$ as a certain solution of the above and the stopping time will be given as usual through

$$T_{\mathbf{AY}}^{\nu,\mu} = \inf\{t \geq 0 : S_t \geq \Psi_\mu^\nu(B_t)\}. \tag{5.6}$$

Let us look first at a simple example. Let $B$ be a Brownian motion starting from zero and $\nu \in \mathcal{M}_0^1$ be concentrated on $a$ and $b$, in other words $\nu \sim B_{T_{a,b}}$. Then $\mu$ is also centered as $U\mu \leq U\nu$. Let $a_\mu = \inf\{t : t \in supp(\mu)\} \leq a$ be the lower bound of the support of $\mu$ and $\Psi_\mu$ be the standard barycentre function defined through (5.2). Writing the explicit formula for $U\nu$ in (5.5) we obtain

$$\Psi_\mu^\nu(x) = \begin{cases} a & , \text{ for } x \leq a_\mu \\ \frac{\int_{[x,\infty)} y d\mu(y) + \frac{ba}{b-a}}{\mu([x,\infty)) + \frac{a}{b-a}} & , \text{ for } a_\mu < x < b \\ \Psi_\mu(x) & , \text{ for } x \geq b. \end{cases}$$

Consider now the general setup. The equation (5.5) has either two solutions $t_x^1 < t_x^2$ or solutions form an interval $t \in [t_x^-, t_x^+]$, $t_x^- \leq t_x^+$, see Figure 1. Let us denote these two situations respectively by $I$ and $II$. The latter happens if and only if the two potentials coincide in $x$. Finally, for a probability measure $\rho$ let $a_\rho = \inf\{t : t \in supp(\rho)\}$, $b_\rho = \sup\{x : x \in supp(\rho)\}$ be respectively the lower and the upper bound of the support of $\rho$. Define:

$$\Psi_\mu^\nu(x) = \begin{cases} a_\nu & , \text{ for } x \leq a_\mu \\ t_x^2 & , \text{ for } a_\mu < x < b_\mu \text{ if } I \\ x & \text{ for } (a_\mu < x < b_\mu \text{ if } II) \text{ and for } (x \geq b_\mu). \end{cases} \tag{5.7}$$

A simple observation is required to see that the definition agrees with the intuition. Since $U\mu \leq U\nu \leq U\delta_m$, and $U\mu$ and $U\nu$ coincide with $U\delta_m$ outside their supports, the support of $\nu$ is contained in the support of $\mu$, $[a_\nu, b_\nu] \subset [a_\mu, b_\mu]$. This grants that $\Psi_\mu^\nu(x) \geq x$.

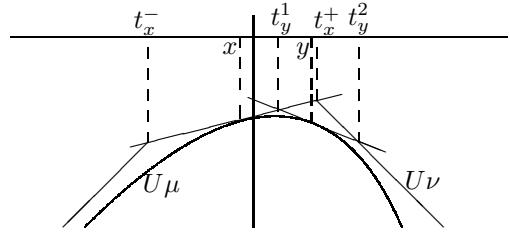

Figure 1: Possible intersections of tangents for nontrivial initial law.

It is also worth mentioning that the function $\Psi_\mu^\nu$ lost some of the properties characterizing $\Psi_\mu$. In particular we can have $\Psi_\mu^\nu(x_1) = x_1$ and $\Psi_\mu^\nu(x_2) > x_2$, for some $x_2 > x_1$.

**Proposition 5.2.** *Let $(M_t : t \geq 0)$ be a continuous, real-valued, local martingale with $\langle M, M \rangle_\infty = \infty$, $M_0 \sim \nu$, and $\mu \in \mathcal{M}^1$ with $U\mu \leq U\nu$. Define $\Psi_\mu^\nu$ and $T_{\mathbf{AY}}^{\nu,\mu}$ through (5.7) and (5.6) respectively. Then $(M_{t \wedge T_{\mathbf{AY}}^{\nu,\mu}})$ is a uniformly integrable martingale and $M_{T_{\mathbf{AY}}^{\nu,\mu}} \sim \mu$.*

The proof of this result is analogous to the proof of Theorem 5.1, taking into account the above discussion of $\Psi_\mu^\nu$.

We remark that the construction presented above might be used to mix or refine some of the known solutions to the embedding problem. For example if $\mu([-a,a]) = 0$ we can first stop at $T_{-a,a}$ and then use $T_{\mathbf{AY}}^{\nu,\mu}$ with $\nu = \frac{1}{2}(\delta_{-a} + \delta_a)$.

### 5.3. Minimizing the minimum - the reversed AY solution

In Section 5.1 we saw how to derive Azéma and Yor's solution to the Skorokhod problem using potential theory. It is easy to see that the same arguments allow us to recover another solution, which we call the reversed Azéma-Yor solution. It was first studied by van der Vecht [121] (see also Meilijson [74]). Recall that $l_x$ was the tangent to the potential $U\mu$ in $(x, U\mu(x))$. We considered its intersection with the line $\{(t, -t) : t \in \mathbb{R}\}$, which was denoted $(\Psi_\mu(x), -\Psi_\mu(x))$. Consider now the intersection of $l_x$ with the line $\{(t, t) : t \in \mathbb{R}\}$ and denote it by $(\Theta_\mu(x), \Theta_\mu(x))$. A simple calculation (analogous to the derivation of (5.2) above) shows that:

$$\Theta_\mu(x) = \frac{1}{\mu((-\infty, x))} \int_{(-\infty, x)} y \, d\mu(y). \tag{5.8}$$

This is a valid formula for $x$ which is not an atom. Looking at the discrete case we see that we want $\Theta_\mu$ to be right continuous, so we finally need to define:

$$\Theta_\mu(x) = \frac{1}{\mu((-\infty, x])} \int_{(-\infty, x]} y \, d\mu(y). \tag{5.9}$$

Denote $s_t = -\inf\{M_s : s \leq t\}$ the infimum process and define the following stopping time

$$T_{\mathbf{rAY}} = \inf\{t \geq 0 : -s_t \leq \Theta_\mu(M_t)\}. \tag{5.10}$$

**Corollary 5.3.** *Let $(M_t)_{t\geq 0}$ be a continuous local martingale, $M_0 = 0$ a.s., with $\langle M, M \rangle_\infty = \infty$ a.s. and let $\mu$ be a centered probability measure on $\mathbb{R}$. Define $\Theta_\mu$ through (5.9). Then $T_{\mathbf{rAY}}$ is a stopping time and $M_{t\wedge T_{\mathbf{rAY}}} \sim \mu$. Furthermore, $(M_{t\wedge T_{\mathbf{rAY}}} : t \geq 0)$ is a uniformly integrable martingale. As a consequence, if $\int_\mathbb{R} x^2 d\mu(x) = v < \infty$ then $\mathbb{E} T_{\mathbf{rAY}} = v$.*

We called this a corollary as this solution is still an Azéma-Yor solution. Indeed if $X \sim \mu$ then $T_{\mathbf{rAY}}$ is just $T_{\mathbf{AY}}$ for $\nu$ the law of $-X$ and the continuous local martingale $N_t = -M_t$. This shows also that we could handle non-trivial initial laws for the reversed Azéma-Yor solution. Speaking informally, in the original Azéma-Yor solution one takes tangents from left to right and in the reversed version one takes the tangents from right to left. We will see below that it yields some symmetric optimal properties.

### 5.4. Martingale theory arguments

In the preceeding, we chose to discuss the Azéma-Yor solution using one-dimensional potential theory. This methodology allowed us to understand well what is really being done and what is really being assumed, opening itself easily to some generalizations, as we saw above and as we will see again in Section 9. One can choose however a different methodology, which has its own advantages as well. We sketch now the heart of the proof that goes through martingale theory. These are the ideas that actually led to the original derivation of the solution in [3].

Let $f : [0, \infty) \to \mathbb{R}$ be an integrable function and denote its primitive by $F(x) = \int_0^x f(y)dy$. Then

$$N_t^f = F(S_t) - (S_t - M_t)f(S_t), t \geq 0 \tag{5.11}$$

is a local martingale, where $M$ and $S$ are as in the previous section. This can be easily seen when $f$ is continuously differentiable, since integrating by parts we obtain

$$\begin{aligned} N_t^F &= F(S_t) - \int_0^t f(S_u)dS_u + \int_0^t f(S_u)dM_u - \int_0^t (S_u - M_u)f'(S_u)dS_u \\ &= F(S_t) - F(S_t) + \int_0^t f(S_u)dM_u = \int_0^t f(S_u)dM_u, \end{aligned} \tag{5.12}$$

as $(S_u - M_u)$ is zero on the support of $dS_u$. The result extends to more general functions $f$ through a limit procedure[7]. Given a centered probability measure $\mu$, we want to find $T$, which embeds $\mu$, of the form $T = \inf\{t \geq 0 : S_t \geq$

---

[7] It can also be verified for $f$ an indicator function and then extended through the monotone class theorem. For more on martingales of the form (5.11) see [85]

$\Psi(M_t)\}$. Our goal is therefore to specify the function $\Psi$. Suppose we can apply the optional stopping theorem to the martingale $N^F$ and the stopping time $T$. We obtain

$$\mathbb{E}\, F(S_T) = \mathbb{E}(S_T - M_T)f(S_T). \qquad (5.13)$$

Note that $M_T \sim \mu$ and $S_T = \Psi(M_T)$. We suppose also that $\Psi$ is strictly increasing. If we denote by $\rho$ the law of $S_T$ and $\overline{\rho} = \rho([x, +\infty))$ its tail, the above equation can be rewritten as

$$\int F(x) d\rho(x) = \int (x - \Psi^{-1}(x))f(x) d\rho(x), \text{ integrating by parts}$$

$$\int f(x) \overline{\rho}(x) dx = \int (x - \Psi^{-1}(x))f(x) d\rho(x), \text{ which we rewrite as}$$

$$d\overline{\rho}(x) = -\frac{\overline{\rho}(x)}{x - \Psi^{-1}(x)} dx, \text{ since } f \text{ was a generic function}.$$

We know that $S_T = \Psi(M_T)$, so that $\overline{\mu}(x) = \mathbb{P}(M_T \geq x) = \mathbb{P}(S_T \geq \Psi(x)) = \overline{\rho}(\Psi(x))$. Inserting this in the above yields

$$d\overline{\mu}(x) = d\Big(\overline{\rho}(\Psi(x))\Big) = -\frac{\overline{\rho}(\Psi(x))}{\Psi(x) - x} d\Psi(x) = -\frac{\overline{\mu}(x)}{\Psi(x) - x} d\Psi(x)$$

$$\Psi(x) d\overline{\mu}(x) - x d\overline{\mu}(x) = -\overline{\mu}(x) d\Psi(x)$$

$$\Psi(x) d\overline{\mu}(x) + \overline{\mu}(x) d\Psi(x) = x d\overline{\mu}(x)$$

$$d\Big(\Psi(x)\overline{\mu}(x)\Big) = x d\overline{\mu}(x)$$

$$\Psi(x)\overline{\mu}(x) = \int_{[x,+\infty)} y d\mu(y),$$

and we find the formula (5.2).

### 5.5. Some properties

The Azéma-Yor solution is known to maximize the supremum in the following sense. Let $\mu$, $M$ and $S$ be as in the previous section, and $R$ any stopping time with $M_R \sim \mu$ and $(M_{R \wedge t} : t \geq 0)$ a uniformly integrable martingale. Then for any $z > 0$, $\mathbb{P}(S_R \geq z) \leq \mathbb{P}(S_{T_{\mathbf{AY}}} \geq z)$. This was observed already by Azéma and Yor [2], as a direct consequence of work by Dubins and Gilat [32], which in turn relied on Blackwell and Dubins [8]. We will try to present here a unified and simplified version of the argumentation. The property was later argued in an elegant way by Perkins [90] and it is also clear from our proof in Section 5.1 (as stressed when discussing measures with finite support, see page 349).

For simplicity assume that $\mu$ has a positive density, so that $\overline{\mu}^{-1}$ is a well defined continuous function. Let $f(x) = \mathbf{1}_{[z,+\infty)}(x)$, then $N_t^F = (S_t - z)^+ - (S_t - M_t)\mathbf{1}_{[z,+\infty)}(S_t)$ is a martingale by (5.11). Applying the optional stopping theorem for $(R \wedge n)$ and passing to the limit (using the monotone convergence theorem and the uniform integrability of $(M_{R \wedge t} : t \geq 0)$) we see that $z \mathbb{P}(S_R \geq$ $z) = \mathbb{E} \, M_R \mathbf{1}_{[z,+\infty)}(S_R)$. Denote $\mathbb{P}(S_R \geq z) = p$. It suffices now to note, that

$$\begin{aligned} zp &= \mathbb{E} \, M_R \mathbf{1}_{[z,+\infty)}(S_R) \leq \mathbb{E} \, M_R \mathbf{1}_{[\overline{\mu}^{-1}(p),+\infty)}(M_R) \\ &= \int_{\overline{\mu}^{-1}(p)}^{\infty} x d\mu(x) = p\Psi_\mu\big(\overline{\mu}^{-1}(p)\big), \quad \text{thus} \\ z &\leq \Psi_\mu\big(\overline{\mu}^{-1}(p)\big), \quad \text{which yields} \\ p &\leq \overline{\mu}\big(\Psi_\mu^{-1}(z)\big) = \mathbb{P}\left(M_{T_{\mathbf{AY}}} \geq \Psi_\mu^{-1}(z)\right) = \mathbb{P}(S_{T_{\mathbf{AY}}} \geq z), \end{aligned}$$

which we wanted to prove. Note that the passage between the last two lines is justified by the fact that $\Psi_\mu^{-1}$ is increasing and $\overline{\mu}$ is decreasing.

This optimality property characterizes the Azéma-Yor solution. Furthermore, it characterizes some of its generalizations - for example in the case of Brownian motion with drift, see Peskir [93]. It also implies an optimality property of the reversed solution presented in Section 5.3. Namely, with $R$ as above, we have $\mathbb{P}(s_R \geq z) \leq \mathbb{P}(s_{T_{\mathbf{rAY}}} \geq z)$, for any $z \geq 0$, where $s_t = -\inf_{u \leq t} M_u$.

The second optimality property says that this solution is pointwise the smallest one. More precisely, if for any stopping time $S \leq T_{\mathbf{AY}}$ and $M_S \sim M_{T_{\mathbf{AY}}}$ then $S = T_{\mathbf{AY}}$ a.s. Actually, it is a general phenomena that we will discuss in Section 8, along with some other general properties of stopping times.

## 6. The $H^1$-embedding

Given a solution $T$ to the Skorokhod stopping problem we can ask about the behavior of supremum or infimum of our process up to time $T$. This question has attracted a lot of attention. In particular some solutions were constructed precisely in order to minimize the supremum and maximize the infimum and conversely.

Let us denote by $(X_t)$ the real valued process we are considering and $\mu$ a probability measure on the real line. The following problem:

*Find some necessary and sufficient conditions on $(X, \mu)$ for the existence of a solution $T$ to the Skorokhod embedding problem such that $X_T^* = \sup\{|X_s| : s \leq T\}$ is in $L^p$*,

is called the $H^p$-embedding problem. We will restrain ourselves to the Brownian setup, even though most of the ideas below work for real-valued diffusions on natural scale. We will present a solution of the $H^p$-embedding by Perkins. He actually gives a solution to the Skorokhod embedding problem, which minimizes $X_T^*$. We then proceed to Jacka's solution, which is a counterpart of Perkins, as it maximizes $X_T^*$. We mention that several authors have investigated in detail the behavior of the maximum of a uniformly integrable martingale and we gathered theirs works in Section 10.3.

### 6.1. A minimal $H^1$-embedding

The problem of $H^p$-embedding was also considered by Walsh and Davis [27] but received a complete and accurate description in Perkins' work in 1985 [90]. Using Doob's $L^p$-inequalities one can see easily that for $p > 1$ an $H^p$-embedding exists if and only if $\int |x|^p d\mu(x) < \infty$, and if $p < 1$ every embedding is in fact an $H^p$-embedding. For $p = 1$ the situation is more delicate, and it goes back to an analytical problem solved by Cereteli [15].

**Theorem 6.1 (Cereteli-Davis-Perkins).** *Consider Brownian motion and a centered probability measure $\mu$. The following are equivalent:*

- *There exists an $H^1$-embedding of $\mu$.*

- *The Skorokhod solution, described in Section 3.1, is an $H^1$-embedding of $\mu$.*

- $H(\mu) = \int_0^\infty \frac{1}{\lambda} |\int_\mathbb{R} x \mathbf{1}_{|x| \geq \lambda} d\mu(x)| d\lambda < \infty.$

Perkins [90] associated with every centered probability measure two functions: $\gamma_-$ and $\gamma_+$ (see the original work for explicit formulae). Then if

$$S_t = \sup_{u \leq t} B_u \quad and \quad s_t = -\inf_{u \leq t} B_u \qquad (6.1)$$

are respectively the supremum and infimum processes we can define:

$$\rho = \inf\{t > 0 : B_t \notin (-\gamma_+(S_t), \gamma_-(s_t))\}. \qquad (6.2)$$

If the target distribution has an atom in 0, we just have to stop in 0 with the appropriate probability. This is done through an external randomization. We take an independent random variable $U$, uniformly distributed on $[0, 1]$ and put

$$T_\mathbf{P} = \rho \mathbf{1}_{U > \mu(\{0\})}, \qquad (6.3)$$

where P stands for Perkins. $T_\mathbf{P}$ is a (randomized) stopping time that embeds $\mu$ and $T_\mathbf{P}$ is an $H^1$-embedding if there exists one. Moreover, Perkins defined his stopping time in order to control the law of the supremum and the infimum of the process and his stopping time is optimal.

**Theorem 6.2 (Perkins).** $(B_{t \wedge T_\mathbf{P}} : t \geq 0)$ *is a uniformly integrable martingale and $B_{T_\mathbf{P}} \sim \mu$. $T_\mathbf{P}$ is an $H^1$-embedding if and only if $H(\mu) < \infty$. Furthermore, for any solution $T$ of the Skorokhod embedding problem for $\mu$, for all $\lambda > 0$,*

$$\mathbb{P}(S_{T_\mathbf{P}} \geq \lambda) \leq \mathbb{P}(S_T \geq \lambda), \qquad (6.4)$$
$$\mathbb{P}(s_{T_\mathbf{P}} \geq \lambda) \leq \mathbb{P}(s_T \geq \lambda), \qquad (6.5)$$
$$\mathbb{P}(B^*_{T_\mathbf{P}} \geq \lambda) \leq \mathbb{P}(B^*_T \geq \lambda). \qquad (6.6)$$

Actually this theorem is not hard to prove - the stopping time is constructed in such a way as to satisfy the above. The hard part was to prove that $T_\mathbf{P}$ embeds $\mu$. The construction of $T_\mathbf{P}$ helps Perkins to give explicit formulae for the three probabilities on the left side of (6.4)-(6.6). Using Perkins' framework we could also define a different stopping time in such a way that it would maximize stochastically the supremum $S$. Of course what we end up with is the Azéma-Yor stopping time (see Sections 5.1 and 5.5) and this is yet an another approach which allows us to recover their solution in a natural way.

The work of Perkins has been generalized by Cox and Hobson [23, 24]. Several authors examined the law of supremum (and infimum) of a martingale with a fixed terminal distribution (see Section 10.3). More generally, the link between the Skorokhod embedding problem and optimal stopping theory was studied. We develop this subject in Section 10.

### 6.2. A maximal $H^1$-embedding

We saw in Section 5 that the Azéma-Yor construction maximizes stochastically the distribution of $S_T$, where $B_T \sim \mu$, fixed, in the class of stopping times such that $(B_{t \wedge T})_{t \geq 0}$ is a uniformly integrable martingale. Moreover, we saw also that its reversed version $T_{\mathbf{rAY}}$, minimizes the infimum. How can we combine these two? The idea came from Jacka [58] and allowed him to construct a maximal $H^1$ embedding. We will describe his idea using the Chacon and Walsh methodology, which is quite different from his, yet helps, in our opinion, to see what is really being done.

We will continue to work within the Brownian setting, even though everything works in exactly the same way for continuous local martingales. We want to find a stopping time $T$ such that $B_T \sim \mu$, $(B_{t \wedge T} : t \geq 0)$ is a uniformly integrable martingale and if $R$ is another such stopping time, then $\mathbb{P}(B_R^* \geq z) \leq \mathbb{P}(B_T^* \geq z)$, for $z \geq 0$. Naturally big values of $B^*$ come either from big values of $S$ or $s$. We should therefore wait "as long as possible" to decide which of the two we choose, and then apply the appropriate extremal embedding. To see whether we should use big negative or big positive values we will stop the martingale at the exit time from some interval $[-k, k]$, which we assume to be symmetric. This corresponds to cutting the initial potential with a horizontal line. The biggest value of $k$ which we can take is obtained by requiring that the horizontal line be tangent to the potential $U\mu$.

For sake of simplicity we suppose that the measure $\mu$ has a positive density. We need to determine the point $x^*$ in which the derivative of $U\mu$ becomes zero. Looking at (5.1) we see that $\overline{\mu}(x^*) = \frac{1}{2}$. In other words, we chose the median of $\mu$ for $x^*$. The value of $k^*$ is then given as the intersection of the tangent to the graph of $U\mu$ in $x^*$ with the line $y = -x$, which is precisely $k^* = \Psi_\mu(x^*)$, where $\Psi_\mu$ is given by (5.2).

As suggested, we first stop at the exit time from $[-k^*, k^*]$. Obviously Brownian motion stops either at $k^*$, in which case it will be finally stopped inside $[U\mu(x^*), \infty)$, or at $-k^*$ in which case it will be finally stopped inside $(-\infty, U\mu(x^*)]$. If we exit at the top, we want to maximize the supremum, and so we will use the standard Azéma-Yor embedding: $T_1 = \inf\{t > T_{-k,k} : S_t \geq \Psi_\mu(B_t)\}$. If we exit at the bottom, we need to minimize the infimum, so we will use the reversed Azéma-Yor solution: $T_2 = \inf\{t > T_{-k,k} : -s_t \leq \Theta_\mu(B_t)\}$. We have then the following proposition,

**Proposition 6.3 (Jacka [58]).** *In the above setup*

$$\mathbb{P}\left(B_R^* \geq z\right) \leq \min\left\{1, \mu([\psi_\mu^{-1}(z), +\infty)) + \mu((-\infty, \Theta_\mu^{-1}(z)])\right\}. \tag{6.7}$$

*Furthermore, this inequality is optimal as equality is attained for a stopping time $T_{\mathbf{J}}$ defined through*

$$T_{\mathbf{J}} = T_1 \wedge T_2 = \inf\left\{t \geq T_{-k^*, k^*} : B_t^* \geq \max(\Psi_\mu(B_t), \Theta_\mu(B_t)\right\}. \tag{6.8}$$

If the measure $\mu$ has atoms the above Proposition is still true taking left inverses of $\Psi_\mu$ and $\Theta_\mu$ and extending correctly the definition of $k^*$. More precisely, if we write $m = \inf\{x : \overline{\mu}(x) \geq \frac{1}{2}\}$ for the median of $\mu$, then if $\mu(\{m\}) > 0$ we have a family of tangents to $U\mu$ in $m$. The value $\Psi_\mu(m)$ results from taking the left-most tangent (cf. Section 5.1) and not the horizontal one. We have therefore $k^* > \Psi_\mu(m)$. The difference is given by the slope of the left-most tangent, $(2\overline{\mu}(m) - 1)$, times the increment interval $(\Psi_\mu(m) - m)$. We obtain therefore: $k^* = \Psi_\mu(m) + (2\overline{\mu}(m) - 1)(\Psi_\mu(m) - m)$, which is the quantity given in Lemma 6 in Jacka [58].

## 7. Root's solution

All the solutions described so far have one thing in common: stopping times are defined through some functionals of Brownian motion up to time $T$. Intuitively, if one searches for a stopping time with minimal variance it should have a very simple history-dependence structure, preferably it should just be the first hitting time of some set. Of course this cannot be just a subset of the space $\mathbb{R}$ but rather of the time-space $\mathbb{R}_+ \times \mathbb{R}$. Root [103] showed that indeed for any centered probability measure $\mu$ on $\mathbb{R}$ with finite second moment, there exists a set $R_\mu \subset \mathbb{R}_+ \times \mathbb{R}$, such that the first entrance time of $(t, B_t)$ into $R_\mu$, embeds $\mu$ in Brownian motion $B$. The first hitting times for the space-time process of Brownian motion correspond to réduites of potential measures and we should not be surprised that to investigate Root's solution Rost [104], [105], [106] needed a fair amount of potential theory.

### 7.1. Solution

The key notion is that of a barrier.

**Definition 7.1.** *A subset $R$ of $[0, +\infty] \times [-\infty, +\infty]$ is a barrier if*

1. *$R$ is closed,*
2. *$(+\infty, x) \in R$ for all $x \in [-\infty, +\infty]$,*
3. *$(t, \pm\infty) \in R$ for all $t \in [0, +\infty]$,*
4. *if $(t, x) \in R$ then $(s, x) \in R$ whenever $s > t$.*

**Theorem 7.2 (Root [103]).** *For any probability measure $\mu$ on $\mathbb{R}$ with $\int x^2 d\mu(x) = v < \infty$ and $\int x d\mu(x) = 0$, there exists a barrier $R_\mu$ such that the first hitting time of the barrier*

$$T_{\mathbf{R}}^\mu = T_{R_\mu} = \inf\{t \geq 0 : (t, B_t) \in R_\mu\}, \tag{7.1}$$

*solves the Skorokhod problem for $\mu$, i.e. $B_{T_{R_\mu}} \sim \mu$ and $\mathbb{E} T_{R_\mu} = v$.*

**Remark on notation:** *We write $T_{\mathbf{R}}^\mu$ for Root's stopping time, which solves the Skorokhod embedding problem for the measure $\mu$ and $T_{\mathbf{R}}$ for a general "Root's stopping time" that is for the first hitting time of some barrier $R$. We write $T_{R_\mu}$ for the first hitting time of a particular barrier $R_\mu$.*

The proof is carried out through a limit procedure. It is not hard to see that a barrier $R_\mu$ exists for a measure with finite support $\mu = \sum_{i=1}^n \alpha_i \delta_{x_i}$. $R_\mu$ is just of the form

$$R_\mu = \left( \bigcup_{i=1}^n [b_i, +\infty] \times \{x_i\} \right) \cup \Big( [0, +\infty] \times [-\infty, x_1] \cup [0, +\infty] \times [x_n, +\infty] \Big). \quad (7.2)$$

The points $(b_i)$ are not given explicitly but only their existence is proven. One then needs some technical observations in order to pass to the limit.

Let $H$ denote the closed half plane, $H = [0, +\infty] \times [-\infty, \infty]$. Root [103] defines a metric $r$ on the space $\mathcal{C}$ of all closed subsets of $H$. Map $H$ homeomorphically to a bounded rectangle $F$ by $(t, x) \to (\frac{t}{1+t}, \frac{x}{1+|x|})$. Let $F$ be endowed with the ordinary Euclidean metric $\rho$ and denote by $r$ the induced metric on $H$. For $C, D \in \mathcal{C}$ put

$$r(C, D) = \max\{\sup_{x \in C} r(x, D), \sup_{y \in D} r(y, C)\}. \quad (7.3)$$

Equipped with $r$, $\mathcal{C}$ is a separable, compact metric space and $\mathcal{R}$, the subspace of all barriers is closed in $\mathcal{C}$ and hence compact. Furthermore, this metric allows us to deal with convergence in probability of first hitting times of barriers. More precisely, the application which associates with a barrier $R$ its first hitting time, i.e. $R \to T_R$, is uniformly continuous from $(\mathcal{R}, r)$ into the set of first hitting times of barriers equipped with convergence in probability.

**Lemma 7.3 (Root [103], Loynes [70]).** *Let $R$ be a barrier with corresponding stopping time $T_R$. If $\mathbb{P}(T_R < \infty) = 1$, then for any $\epsilon > 0$ there exists $\delta > 0$ such that if $R_1 \in \mathcal{R}$ and $r(R, R_1) < \delta$ then $\mathbb{P}(|T_R - T_{R_1}| > \epsilon) < \epsilon$, where $T_{R_1}$ is defined via (7.1). If $\mathbb{P}(T_R = \infty) = 1$ then for any $K > 0$ there exists $\delta > 0$ such that $r(R, R_1) < \delta$ implies $\mathbb{P}(T_{R_1} < K) < \epsilon$.*

*If a sequence of barriers converges, $R_n \xrightarrow[n \to \infty]{r} R$ with $\mathbb{E} T_{R_n} < K < \infty$, then $\mathbb{E} T_R \leq K$ and $\mathbb{P}(|T_{R_n} - T_R| > \epsilon) \xrightarrow[n \to \infty]{} 0$.*

With this lemma the theorem is proven taking a sequence $\mu_n$ of probability measures with finite supports converging to $\mu$, $\mu_n \Rightarrow \mu$.

### 7.2. Some properties

We start with some easy properties and then pass to the more involved results of Rost. The first property, established by Root [103] for his stopping time, asserted the equivalence between integrability of $B_{T_{\mathbf{R}}}^{2k}$ and $T_{\mathbf{R}}^k$, for $k \in \mathbb{N}$. We quoted it, in a more general form, in Preliminaries (Proposition 2.1) as it is true for any stopping time $T$ such that $(B_{t \wedge T} : t \geq 0)$ is a uniformly integrable martingale.

Several further properties were described by Loynes [70], one of which is particularly important for us. In order to state it, we need to introduce the notion of a regular barrier. Let $R \in \mathcal{R}$ be a barrier. Define

$$x_+^R = \inf \left\{ y \in [0, \infty] : [0, +\infty] \times \{y\} \subset R \right\}$$
$$x_-^R = \sup\{y \in [-\infty, 0] : [0, +\infty] \times \{y\} \subset R\}, \tag{7.4}$$

the first coordinates above and below zero respectively at which $R$ is a horizontal line.

**Definition 7.4.** *We say that a barrier $R \in \mathcal{R}$ is regular if $[0, +\infty] \times [x_+^R, +\infty] \subset R$ and $[0, +\infty] \times [-\infty, x_-^R] \subset R$.*
*We say that two barriers $R_1, R_2$ are equivalent, $R_1 \sim R_2$, if $x_-^{R_1} = x_-^{R_2}$, $x_+^{R_1} = x_+^{R_2}$ and the two barriers coincide on $[0, +\infty] \times [x_-^{R_1}, x_+^{R_1}]$.*

Thus a regular barrier $R$ is one that is composed of horizontal lines above $x_+^R$ and below $x_-^R$. Clearly each barrier is equivalent to exactly one regular barrier. Furthermore, for a barrier $R$ we have that $T_R \le T_{x_-^R, x_+^R}$ and so stopping times corresponding to two equivalent barriers $R_1 \sim R_2$ are equal $T_{R_1} = T_{R_2}$. This shows that we can limit ourselves to regular barriers.

**Proposition 7.5 (Loynes [70]).** *Let $\mu$ be a centered probability measure on $\mathbb{R}$ with $\int x^2 d\mu(x) = v < \infty$. Then there exists exactly one regular barrier $R_\mu$, such that $B_{T_{R_\mu}} \sim \mu$ and $\mathbb{E}\, T_{R_\mu} = v$.*

Loynes showed also that a wider class of probability measures $\mu$ can be embedded in Brownian motion through Root's construction. However, the uniqueness result does not carry over to this more general setup.

The fundamental property of Root's construction is that it minimizes the variance of the stopping time. This conjecture, made by Kiefer [64], was proved by Rost.

**Theorem 7.6 (Rost [107]).** *Let $\mu$ be a centered probability measure on $\mathbb{R}$ with $\int x^2 d\mu(x) = v < \infty$ and let $T_{\mathbf{R}}^\mu$ denote Root's stopping time embedding $\mu$. Assume that $\mathbb{E}(T_{\mathbf{R}}^\mu)^2 < \infty$. Then, for any stopping time $S$ such that $B_S \sim \mu$ and $\mathbb{E}\, S = v$, we have: $\mathbb{E}(T_{\mathbf{R}}^\mu)^2 \le \mathbb{E}\, S^2$ and $\mathbb{E}(T_{\mathbf{R}}^\mu)^2 = \mathbb{E}\, S^2$ if and only if $S = T_{\mathbf{R}}^\mu$ a.s.*
*More generally, for any $0 < p < 1$, $\mathbb{E}(T_{\mathbf{R}}^\mu)^p \ge \mathbb{E}\, S^p$ and for any $p > 1$, $\mathbb{E}(T_{\mathbf{R}}^\mu)^p \le \mathbb{E}\, S^p$.*

We will try to explain briefly, and at an intuitive level, how to arrive at the above result. Actually, Rost considers the more general setup of a Markov process (that is assumed transient for technical reasons but the extension is fairly simple).

As in Section 2.2, let $(X_t : t \ge 0)$ be a Markov process (relative to its natural filtration $(\mathcal{F}_t)$) on a compact metric space $E$ associated with a semigroup $(P_t^X)$ and a potential kernel $U^X$. Let $\nu$ and $\mu$ be two probability measures on the Borel sets of $E$ with $\sigma$-finite potentials $\mu U^X$ and $\nu U^X$. The measure $\nu$ is the initial law: $X_0 \sim \nu$ and $\mu$ is the law we want to embed. Recall the characterization of the balayage order given in Section 3.5. For existence of an embedding, the measures have to satisfy $\mu U^X \leq \nu U^X$. Then, there exists a stopping time $T$ which embeds $\mu$, $\mu = \nu P_T^X$.

We are looking for a stopping time $T$ with minimal variance. It is shown actually that it satisfies the stronger property of minimizing $\mathbb{E}_\nu \int_{t \wedge T}^T f(X_s) ds$ for Borel, positive functions $f$, and any $t > 0$. This in turn is equivalent to minimizing $\nu P_{t \wedge T}^X U^X$ for all $t > 0$. The idea is to look at the space-time process. We consider a measure $N$ on $E \times \mathbb{R}$ defined through

$$N(\mathcal{A} \times \mathcal{B}) = \int_\mathcal{B} N_t(\mathcal{A}) dt, \text{ where } N_t = \nu U^X \mathbf{1}_{t<0} + \nu P_{t \wedge T}^X U^X \mathbf{1}_{t \geq 0}.$$

We take its réduite $\hat{N}$ (with respect to the semigroup of the space-time process) and show that it is of the form

$$\hat{N}(\mathcal{A} \times \mathcal{B}) = \int_\mathcal{B} \hat{N}_t(\mathcal{A}) dt,$$

with $\hat{N}_{\cdot}$ decreasing, right continuous and $\hat{N}_t = \nu P_{T_R \wedge t}^X U^X$ for a certain stopping time with $\nu P_{T_R}^X = \mu$. It follows that $T_R$ minimizes $\nu P_{t \wedge T}^X U^X$. The last step consists in showing that $T_R$ is actually Root's stopping time for a certain barrier $R$. This is not surprising after all, as we saw that réduites are realized as hitting times of some sets and Root's stopping time is just a hitting time for the space-time process.

### 7.3. The reversed barrier solution of Rost

In this section we describe a solution proposed by Rost, which is a counterpart for Root's solution as it maximizes the variance of the stopping time. We have not found it published anywhere and we cite it following the notes by Meilijson [74]. The basic idea is to look at reversed barriers. Intuitively we just substitute the condition that if a point belongs to a barrier then all the points to the right of it also belong to the barrier, by a condition that all the points to the left of it also belong to the barrier.

**Definition 7.7.** *A subset $\rho$ of $[0, +\infty] \times [-\infty, +\infty]$ is called a reversed barrier if*

1. *$\rho$ is closed,*
2. *$(0, x) \in \rho$ for all $x \in [-\infty, +\infty]$,*
3. *$(t, \pm\infty) \in \rho$ for all $t \in [0, +\infty]$,*
4. *if $(t, x) \in \rho$ then $(s, x) \in \rho$ whenever $s < t$.*

**Theorem 7.8 (Rost).** *For any probability measure $\mu$ on $\mathbb{R}$ with $\mu(\{0\}) = 0$, $\int x^2 d\mu(x) = v < \infty$ and $\int x d\mu(x) = 0$, there exists a reversed barrier $\rho_\mu$ such that the stopping time*

$$T_{\rho_\mu} = \inf\{t > 0 : (t, B_t) \in \rho_\mu\}, \tag{7.5}$$

*solves the Skorokhod problem for $\mu$, i.e. $B_{\rho_\mu} \sim \mu$ and $\mathbb{E} T_{\rho_\mu} = v$. Furthermore, if $S$ is a stopping time such that $B_S \sim \mu$ and $\mathbb{E} S = v$ then for any $p < 1$ $\mathbb{E} T_{\rho_\mu}^p \leq \mathbb{E} S^p$, and for any $p > 1$ $\mathbb{E} T_{\rho_\mu}^p \geq \mathbb{E} S^p$.*

Note that there is a certain monotonicity in the construction of a reversed barrier, which is not found in the construction of a barrier. For a point $x$ in the support of $\mu$, say $x > 0$, we must have the following property: there exists some $t_x$ such that $(t_x, x) \in \rho_\mu$ and $(t_x, y) \notin \rho_\mu$ for any $y \in [0, x)$. In other words, the reversed barrier $\rho_\mu$, on the level $x$, has to extend "more to the right" than anywhere below, between 0 and $x$. This, in a natural way, could lead us to revise the notion of a reversed barrier and introduce regular reversed barriers, which would satisfy $\forall x > 0$ if $(t, x) \in \rho$ then $(t, y) \in \rho$, $\forall y \in [x, +\infty)$ and likewise for $x < 0$. We could then obtain uniqueness results for Rost's embedding in an analogy to Proposition 7.5 above. We do not pursue these matters here, but we note that, according to our best knowledge, they have never been investigated.

### 7.4. Examples

Explicit Root's or Rost's solutions are very rarely known, still we can give some examples. First, a trivial example deals with Gaussian measures. If $\mu = \mathcal{N}(0, v)$ then the barrier $R_v = \{(t, x) : t \geq v\}$ yields a deterministic stopping time $T_{R_v} = v$, which naturally embeds the measure $\mu$.

More sophisticated examples deal with square root boundaries. Shepp [113] and Breiman [10], independently, considered the following stopping times[8]:

$$\sigma_{a,c} = \inf\{t \geq 0 : |B_t| = c\sqrt{t+a}\}.$$

These are easily seen to be Rost's stopping times as $\sigma_{a,c} = T_{\rho_{a,c}}$, where $\rho_{a,c} = \{(t, x) \in \mathbb{R}_+ \times \mathbb{R} : x^2 \geq c(t+a)\}$ is a reversed barrier. As $B_{\sigma_{a,c}} = c\sqrt{\sigma_{a,c} + a}$, to characterize the law of $B_{\sigma_{a,c}}$, it suffices to characterize the law of $\sigma_{a,c}$, which in turn is equivalent to giving its Laplace transform. This can be done using exponential martingales or through a reduction to the first hitting times of a given level for a radial Ornstein-Uhlenbeck process. We do not give the explicit formulae and refer to Göing-Jaeschke and Yor [45] (see for example the formulae on page 322) for a recent and extensive study which includes the subject.

Novikov [82] considered the following stopping time:

$$\tau_{a,b,c} = \inf\{t \geq 0 : B_t \leq -a + b\sqrt{t+c}\}, \tag{7.6}$$

---

[8]More precisely, Breiman [10] worked with $T = \inf\{t \geq 1 : B_t \geq c\sqrt{t}\}$ conditionally on $\{B_1 = 0\}$.

where $a \geq 0, c \geq 0$ and $bc^{1/2} < a$. This is Root's stopping time for a barrier $R_{a,b,c} = \{(t,x) \in \mathbb{R}_+ \times \mathbb{R} : x \geq -a + b\sqrt{t+c}\}$, $T_{R_{a,b,c}} = \tau_{a,b,c}$. Again, characterizing the law of $B_{\tau_{a,b,c}}$ is equivalent to characterizing the law of $\tau_{a,b,c}$ itself. To do so Novikov used a general family of exponential martingales of Brownian motion and the optional stopping theorem. In particular Novikov showed that $\mathbb{E}\,\tau_{a,b,c}^v < \infty$ if and only if $b > z_{2v}$, where $z_v$ is the largest root of the equation $D_v(z) = 0$, for $D_v$ the parabolic cylinder function of the parameter $v$.

The stopping times $\sigma_{a,c}$ and $\tau_{a,b,c}$ served Davis [26] to attain the optimal lower and upper bounds on the proportion $\mathbb{E}\,|B_T|^p/\mathbb{E}\,T^{p/2}$, where $T$ is a stopping time, $\mathbb{E}\,T^{p/2} < \infty$, $p > 1$. This is not surprising, as we would expect that to attain the bounds one would have to consider stopping times which, with distribution of $B_T$ fixed, minimize and maximize the expectation $\mathbb{E}\,T^{p/2}$, that is Root's and Rost's stopping times (cf. Theorems 7.6 and 7.8).

Finally we mention two more papers which involve explicit calculations for Root's and Rost's stopping times. For a barrier $R$ it is naturally possible that $\mathbb{P}(T_R = \infty) > 0$. Robbins and Sigmund [100] calculate the probabilities $\mathbb{P}(T_R = \infty)$ or $\mathbb{P}(T_R \geq u)$, $u$ fixed, for the barriers of the type $R = \{(t,x) : x \geq f(t)\}$, where $f$ is a function from a special family, which includes functions such as $f(t) = ct + a$ or $f(t) = \sqrt{(t+\epsilon)\log(t+\epsilon)}$. Lebedev [69] has calculations, for an arbitrary $d$-dimensional diffusion, of the first moment of Rost's stopping times which are the first passage times through a time-dependent, smooth boundaries.

### 7.5. Some generalizations

The problem we want to solve is a certain generalization of the one solved by Rost. In [107] Rost showed that among all stopping times which solve a given Skorokhod embedding problem, the stopping time defined by Root has minimal second moment (which is equivalent to having minimal variance). In the Skorokhod embedding problem, the distribution of $B_T$ is fixed. Here we only fix two moments of this distribution and want to find the stopping time with the minimal variance. That is:

**Problem.** *Let us fix $4 > p > 1$ and a constant $c_p > 0$. Among all stopping times $T$ such that $\mathbb{E}B_T = 0$, $\mathbb{E}T = \mathbb{E}B_T^2 = 1$ and $\mathbb{E}|B_T|^p = c_p$, find the element $T_{min}$ with minimal variance.*

This problem is also linked with the optimal bounds on $\mathbb{E}\,|B_T|^p/\mathbb{E}\,T^{p/2}$ obtained by Davis [26], which we mentioned in the previous section. Davis [26] gives a solution to our problem without the restriction $\mathbb{E}\,|B_T|^p = c_p$.

Unfortunately, we did not succeed in solving this problem. We know it has a solution and we will present a certain conjecture about it.

**Proposition 7.9.** *The stopping time $T_{min}$ exists.*

*Proof.* Recall Definition 7.1 and (7.1). It is an immediate observation that the stopping time we are looking for will be Root's stopping time $T_R$ for some barrier $R$. Indeed, if the $\mu \sim B_{T_{min}}$ is the distribution of $B_{T_{min}}$, then there exists a barrier $R_\mu$ such that Root's stopping time $T_{R_\mu}$ yields $B_{T_{R_\mu}} \sim \mu$. Moreover, from Theorem 7.6, we know that $\mathbb{E}T_{R_\mu}^2 \leq \mathbb{E}T_{min}^2$ and so $\mathbb{E}T_{R_\mu}^2 = \mathbb{E}T_{min}^2$ from the definition of $T_{min}$. Therefore, using Theorem 7.6 again, we conclude that $T_{R_\mu} = T_{min}$ a.s.

We can limit ourselves to Root's stopping times. Let

$$\mathcal{A}_{c_p} = \left\{ \mu \in M_1(\mathbb{R}) : \int x\mu(dx) = 0, \int x^2 \mu(dx) = 1, \int |x|^p \mu(dx) = c_p \right\}.$$

For every $\mu \in \mathcal{A}_{c_p}$ there exists a unique regular barrier $R_\mu$ such that $B_{T_{R_\mu}} \sim \mu$. Let $\mathcal{R}_{c_p}$ be the subset of regular barriers corresponding to the measures in $\mathcal{A}_{c_p}$.

Let $H$ denote the closed half plane and recall the definition (7.3) of the metric $r$, under which $\mathcal{C}$ the space of all closed subsets of $H$, is a separable, compact metric space and $\mathcal{R}$, the subspace of all barriers, is closed in $\mathcal{C}$ and hence compact.

We can introduce an order relation in $\mathcal{R}_{c_p}$ given by: $C \preccurlyeq D \Leftrightarrow \mathbb{E}T_C^2 \leq \mathbb{E}T_D^2$, for $C, D \in \mathcal{R}_{c_p}$. If we prove that there exists $R_{min}$ - a minimal element of $\mathcal{R}_{c_p}$ relative to $\preccurlyeq$ - then, we will have $T_{min} = T_{R_{min}}$. In other words

*If there exists a $R_{min} \in \mathcal{R}_{c_p}$ such that*

$$\mathbb{E}T_{R_{min}}^2 = \inf_{R \in \mathcal{R}_{c_p}} \mathbb{E}T_R^2, \qquad \text{then } T_{min} = T_{R_{min}}.$$

To complete the proof we need to show that $\mathcal{R}_{c_p}$ is closed under $r$. Indeed, choose a sequence of barriers $R_n \in \mathcal{R}_{c_p}$, such that $\mathbb{E}T_{R_n}^2 \searrow \inf_{R \in \mathcal{R}_{c_p}} \mathbb{E}T_R^2$. We can suppose that $R_n$ converges (possibly choosing a subsequence) in $r$ to some barrier $R$. To conclude that $R = R_{min}$ we need to know that $R \in \mathcal{R}_{c_p}$, which follows if $\mathcal{R}_{c_p}$ is a closed subset of $\mathcal{R}$.

In order to prove that $\mathcal{R}_{c_p}$ is a closed subset of $\mathcal{R}$ we essentially follow the ideas used by Root in the proof of Theorem 7.2. Take any sequence $R_n \xrightarrow{r} R$ with $R_n \in \mathcal{R}_{c_p}$. Then $\mathbb{E}T_{R_n} = 1$ for all $n$. It follows from Lemma 7.3 that for any $\epsilon > 0$, $\mathbb{P}(|T_R - T_{R_n}| > \epsilon) \xrightarrow[n \to \infty]{} 0$ and $\mathbb{E}T_R < \infty$. Choosing a subsequence, we may assume that $T_{R_n} \to T_R$ a.s. Furthermore, we can always construct a measure $\mu \in \mathcal{A}_{c_p}$ with $\int x^4 d\mu(x) < \infty$, which implies $\mathbb{E}T_{R_\mu}^2 < \infty$ by Proposition 2.1, so that we may assume the sequence $T_{R_n}$ is bounded in $L^2$, and therefore uniformly integrable, hence converges in $L^1$. In particular $ET_\mathbf{R} = 1$. Finally, Proposition 2.1 also grants the existence of a universal constant $C_4$, such that $\mathbb{E}B_T^4 \leq C_4 \mathbb{E}T^2$, which implies that the sequence $(B_{T_{R_n}})$ is bounded in $L^4$. It converges almost surely, thanks to the continuity of paths and therefore it converges in $L^p$. We see then, that $\mathbb{E}|B_{T_R}|^p = c_p$. This implies that $R \in \mathcal{R}_{c_p}$ and proves that $\mathcal{R}_{c_p}$ is a closed set. We conclude that the stopping time $T_{min}$ exists and is equal a.s. to $T_R$ for a certain barrier $R$. □

We do not know how to obtain $T_{min}$ explicitly. We do not even know if it is unique. Still, we want to give a simple observation, which motivates the conjecture below. We will let $c_p$ change now, so we will write $T_{min} = T_{min}^{c_p} = T_{R_{c_p}}$ and $v_{c_p} := \mathbb{E}T_{R_{c_p}}^2$.

Put $c_p^* = \mathbb{E}|B_1|^p$, then the minimal stopping time is trivial $T_{min}^{c_p^*} = 1$ and it is unique. It corresponds to the vertical barrier $N = \{(x,t) : t \geq 1\}$. Take any $c_p > 0$ and define an order relation in the set $\mathcal{R}_{c_p}$ by: $C \prec D$ iff $r(C, N) \leq r(D, N)$. The minimal element $R_{min}^N$ in $\mathcal{R}_{c_p}$, relative to this order, exists as $\mathcal{R}_{c_p}$ is closed. Denote the stopping time corresponding to this barrier by $T_{min}^N = T_{R_{min}^N}$. The above observation suggested us the following conjecture:

**Conjecture:** *The two minimal stopping times coincide a.s.: $T_{min}^N = T_{min}$.*

This conjecture yields a manageable description of the stopping time - solution to our problem. This comes from the fact that it can be translated somehow into information about the distribution of $B_{T_{min}^{c_p}}$. However, we do not know how to prove it. We only have some ideas that we present now.

Assume that $1 < p < 2$. Consider $f : [0,1) \to [1, \infty)$, $f(c_p) = v_{c_p}$. We will now show that $f$ is a continuous function. We start with the right-continuity. Fix $c_p < 1$ and take a sequence $c_p^n$ decreasing to $c_p$. We have to show that for each subsequence $c_p^{n_k}$ there exists a sub-subsequence $c_p^{n_{k_l}}$ such that $v_{c_p^{n_{k_l}}}$ converges to $v_{c_p}$. Take an arbitrary subsequence of $c_p^n$ and denote it still $c_p^n$. For simplicity of notation put $c := c_p$, $v := v(c_p)$, $R = R_{c_p}$, $c_n = c_p^n$ and $v_n = v_{c_p^n}$. It is an easy observation that $v_n$ cannot converge to anything smaller than $v$ as such a situation would yield a stopping time with smaller variance than $T_{min}^{c_p}$. So it is enough to find a subsequence of barriers $R_{n_k}$ converging to $R$, $R_{n_k} \in \mathcal{R}_{c_{n_k}}$, such that the corresponding stopping times $T_{R_{n_k}}$ will converge in $L^2$.

Let $\mu$ be the distribution of $B_{T_{min}^c}$. Define a family of probability measures by

$$\mu_s = \frac{s}{2}(\delta_{\{-1\}} + \delta_{\{1\}}) + (1-s)\mu, \text{ for } 0 < s < 1.$$

Observe that $\int x\mu_s(dx) = 0$, $\int x^2\mu_s(dx) = 1$ and $\int |x|^p\mu_s(dx) = c + s(1-c) > c$. Corresponding to our sequence of $(c_n)$ we have a decreasing sequence $(s_n)$, $s_n \searrow 0$, $\int |x|^p\mu_{s_n}(dx) = c_n$. Let $R_n$ be the sequence of barriers such that $B_{T_{R_n}} \sim \mu_{s_n}$. Note that $R_n \in \mathcal{R}_{c_n}$. We need also to point out that $\int x^4\mu_s(dx) \leq \int x^4\mu(dx) < \infty$ (this implies that the family $T_{R_n}$ is uniformly bounded in $L^2$). Now we can choose a converging subsequence $R_{n_k}$ and the corresponding stopping times converge with $\mathbb{E}T_{R_{n_k}}^2 \to v$.

To prove the left-continuity of $f$ we follow the same reasoning, we just have to define the appropriate measures which will play the role of $\mu_s$. Fix a constant $a > 1$ so that $\frac{1}{a^{2-p}} < c$. Define a measure $\nu = (1 - \frac{1}{a^2})\delta_{\{0\}} + \frac{1}{2a^2}(\delta_{\{-a\}} + \delta_{\{a\}})$, which plays the role of $\frac{1}{2}(\delta_{\{-1\}} + \delta_{\{1\}})$ above. We defined the measure so as to have: $\int x\nu(dx) = 0$, $\int x^2\nu(dx) = 1$ and $\int |x|^p\nu(dx) = \frac{1}{a^{2-p}}$. Now, for $0 < s < 1$, we take $\eta_s = s\nu + (1-s)\mu$, where $\mu \sim B_{T_{min}^c}$. It follows that $\int x\eta_s(dx) = 0$, $\int x^2\eta_s(dx) = 1$ and $\int |x|^p\eta_s(dx) = c + s(\frac{1}{a^{2-p}} - c) < c$.

It is an easy observation that $f(c) \geq 1$ and $f(c) = 1$ if and only if the variance of the corresponding stopping time $T_{min}^c$ is 0, that is if and only if $T_{min}^c = t$ a.s. for some $t > 0$. Then, since $\mathbb{E}T_{min}^c = 1$, we have $T_{min}^c = 1$ a.s. and so $c = \mathbb{E}|B_1|^p$. In other words, $f$ has one global minimum. We believe that $f$ is a convex function.

## 8. Shifted expectations and uniform integrability

So far, throughout this survey, we have considered the case when $U\mu \leq U\nu$, where $\nu$ is the initial law and $\mu$ is the measure we want to embed (recall Definition 2.2 which gives $U\mu$). A natural question to ask is: what about other probability measures? Can they be embedded into a continuous local martingale $(M_t : t \geq 0)$ as well? The answer is yes. But we have to pay for this - the stopped martingale will no longer be uniformly integrable.

If $(M_{t\wedge T} : t \geq 0)$ is uniformly integrable, $T < \infty$ a.s., then $\mathbb{E} M_0 = \mathbb{E} M_T$. Suppose both $\nu$ and $\mu$ have finite expectations but $\int x d\nu(x) \neq \int x d\mu(x)$. This corresponds to the simplest situation when we do not have the inequality between the two potentials: $U\mu \not\leq U\nu$. Then, if $T$ embeds $\mu$, the martingale $(M_{t\wedge T} : t \geq 0)$ cannot be uniformly integrable, as $\mathbb{E} M_0 \neq \mathbb{E} M_T$. The same applies to the situation when $\int |x| d\nu(x) \leq \infty$ and $\int |x| d\mu(x) = \infty$.

It is very easy however to point out a number of possible stopping rules which embed $\mu$, in the case of $\int |x| d\mu(x) < \infty$. It suffices to note that our usual assumption: $\langle M, M\rangle_\infty = \infty$ a.s., implies that $T_a = \inf\{t \geq 0 : M_t = a\}$ is finite a.s. for any $a \in \mathbb{R}$. We can therefore first wait till $M$ hits $m = \int x d\mu(x)$ and then use any of the embeddings described so far applied to the shifted process starting at $m$.

We have to see therefore two more properties. Firstly, that any probability measure with $\int |x| d\mu(x) = \infty$ can be embedded in $M$, and secondly that in case of $\int |x| d\mu(x) < \infty$ but $U\mu \not\leq U\nu$ we necessarily loose the uniform integrability. Resolution of our first problem is trivial and was proposed by Doob as noted on page 331 (see also exercise VI.5.7 in Revuz and Yor [99]). A different solution is found in Cox and Hobson [23] who generalized the Perkins stopping rule in such a way that it works with any distribution on $\mathbb{R}$. The construction itself is rather technical and we will not present it here. We refer to Section 3.20 and to the original paper for details.

The second question asked above is also simple[9]. Suppose that $M_0 \sim \nu$ and $M_T \sim \mu$. Then for any $x \in \mathbb{R}$ the process $(M_{t\wedge T} - x : t \geq 0)$ is a martingale and therefore $(|M_{t\wedge T} - x| : t \geq 0)$ is a submartingale. In particular the expectations of the latter increase and therefore $-U\nu(x) = \mathbb{E}|M_0 - x| \leq \mathbb{E}|M_T - x| = -U\mu(x)$ for any $x \in \mathbb{R}$, or equivalently $U\mu \leq U\nu$. We have thus proved the following proposition:

**Proposition 8.1.** *Let $(M_t : t \geq 0)$ be a continuous, local, real-valued martingale with $\langle M, M\rangle_\infty = \infty$ a.s. and $M_0 \sim \nu$. For any probability measure $\mu$ on $\mathbb{R}$ there exists a stopping time in the natural filtration of $M$ which embeds $\mu$, that is $M_T \sim \mu$. Furthermore, the stopping time can be taken so that $(M_{t\wedge T} : t \geq 0)$ is a uniformly integrable martingale if and only if $U\mu \leq U\nu$.*

---

[9]Strangely enough we found it treated only in notes by Meilijson [74] sent to us by the author together with remarks on the preliminary version of this work.

A natural question therefore is how to express the idea that the stopping time should be "small", if we cannot use the uniform integrability condition? This brings us to the question of various types, or properties, of stopping times. We gather and discuss them now.

**Definition 8.2.** *A stopping time $T$ is called **minimal** if for any stopping time $S \leq T$, $M_S \sim M_T$ implies $S = T$ a.s. (cf. Monroe [78], Cox and Hobson [24]). A stopping time $T$ is called **standard** if there exists a sequence of bounded stopping times $T_n$, with $\lim T_n = T$ a.s. and $\lim U\mathcal{L}(M_{T_n})(x) = U\mathcal{L}(M_T)(x) > -\infty$ for all $x \in \mathbb{R}$ (cf. Chacon [17], Chacon and Ghoussoub [18], Falkner [36]). A stopping time $T$ is called **ultimate** if for any probability measure $\rho$ with $U\rho \geq U\mu$, where $M_T \sim \mu$, there exists a stopping time $S \leq T$ such that $M_S \sim \rho$ (cf. Meilijson [72], van der Vecht [121]).*

For $T$ a standard stopping time we always have $\mathbb{E}\, M_T = 0$ (as $-\infty < U\mathcal{L}(M_T) < U\delta_0$), while minimal stopping times exist in a broader context and they seem fundamental for the extension of the Skorokhod embedding to non-centered cases. The following proposition explains that the minimality concept also applies in the centered case:

**Proposition 8.3.** *Let $T$ be a stopping time such that $\mathbb{E}\, M_T = 0$. The following are equivalent*

1. *$T$ is minimal,*

2. *$T$ is standard,*

3. *$(M_{t \wedge T} : t \geq 0)$ is a uniformly integrable martingale.*

The equivalence between 1. and 3. was first obtained by Monroe [78], who made an extensive use of the theory of barriers (see Section 7). It was then argued in a much simpler way by Chacon and Ghoussoub [18]. They showed, using one dimensional potential theory, that a minimal stopping time, for centered final distribution, is the same as a standard stopping time (this was also shown by Falkner [36]). Therefore one can deduce properties of minimal stopping times from those of standard stopping times, which in turn are easy to establish (see Chacon [17]).

An immediate conclusion from the above Proposition, is that any solution to the Skorokhod embedding problem in the classical setup is minimal. It seems that in the non-centered setup minimality is the right condition to impose on the stopping times (Cox and Hobson [24]). The following proposition assures that it is also a plausible condition:

**Proposition 8.4 (Monroe [78]).** *For any stopping time $S$ there exists a minimal stopping time $T \leq S$, with $M_S \sim M_T$.*

We end this section with a result about ultimate stopping times. Namely, we treat the question: is it possible for a stopping time to be standard and ultimate at the same time? A negative answer is given by the following

**Proposition 8.5 (Meilijson [72], van der Vecht [121]).** *If $T$ is a standard and ultimate stopping time and $\mu \sim B_T$ then $\mu(\{a,b\}) = 1$, for some $a < 0 < b$.*

The proof consists in showing that such a solution would have to coincide with both the standard Azéma and Yor solution and its reversed version (see Section 5), which in turn are shown to coincide only in trivial cases.

## 9. Working with real-valued diffusions

As the reader may have have observed, the Chacon and Walsh methodology does not rely on any special property of Brownian motion. We underlined this in Section 3.8, and in Section 5.1 while proving the Azéma-Yor construction. There is only one crucial property that was used: the way potential changes if we stop our process at the exit time of $[a, b]$, which is described in *(vii)* of Proposition 2.3. This is also the only property used by Zaremba [126] in his proof of the Azéma-Yor solution.

This property stays true however, for any real-valued diffusion on natural scale. Therefore, generalizing this framework to real-valued diffusions with continuous, strictly increasing scale functions is easy. Hall [50] (see Section 3.4) treated Brownian motion with drift. Recurrent diffusions were treated in the original work of Azéma and Yor [3], and the general case was considered by Pedersen and Peskir [89]. We present here a complete description following mainly a recent paper by Cox and Hobson [23]. We note that more general, not necessarily continuous processes, were treated as well. In Section 3.15 we presented the solution developed by Bertoin and Le Jan [6], and in Section 3.21 we discussed the solution for the age process of Brownian excursions developed by Obłój and Yor [84].

Let $(X_t)_{t \geq 0}$ be a regular (time-homogeneous) diffusion taking values in an interval $I \subset \mathbb{R}$, $X_0 = 0$ a.s. Then, there exists a continuous, strictly increasing scale function $s : I \to \mathbb{R}$ and $M_t = s(X_t)$ is a local martingale on $s(I)$. Furthermore, we can choose $s(0) = 0$ (see for example Revuz and Yor [99] p. 301).

We transform the embedding problem for $X$ into an embedding problem for $M$. Let $\mu$ be a probability measure on $I$, which we want to embed in $X$. Define a measure $\mu^s$ through

$$\mu^s(B) = \mu(s^{-1}(B)), \quad B \in \mathcal{B}(s(I)),$$

that is if a random variable $\xi \sim \mu$ then $s(\xi) \sim \mu^s$. It is easily seen that if we find a stopping time $T$ such that $M_T \sim \mu^s$ then $X_T \sim \mu$ and the embedding problem for $X$ will be solved.

At first glance one would like to conclude applying our previous results to $M$. However some care is needed as $s(I)$ might be bounded and the local martingale $M$ will not be recurrent.

We assume that $\int |x| d\mu^s(x) = \int |s(u)| d\mu(u) < \infty$ and write $m = \int x d\mu^s(x) = \int s(u) d\mu(u)$. There are three distinct cases: $s(I)^\circ = \mathbb{R}$, $s(I)^\circ = (-\infty, \alpha)$ or $s(I)^\circ = (\alpha, \infty)$ and $s(I)^\circ = (\alpha, \beta)$, which we now discuss.

1. $\mathbf{s(I)}^\circ = \mathbb{R}$. This is a trivial situation since the local martingale $(M_t : t \geq 0)$ is recurrent and we can use any standard construction, such as the Azéma-Yor solution (see Section 5), to embed $\mu^s$ in $M$.

2. $\mathbf{s(I)}^\circ = (-\infty, \alpha)$. In this case $(M_t : t \geq 0)$ is a local martingale that is bounded from above by $\alpha$. Suppose $T$ is a stopping time such that $M_T \sim \mu^s$ and let $\gamma_n$ by a localizing sequence for $M$. Then, by Fatou lemma,
$$(-m) = \mathbb{E}(-M_T) = \mathbb{E}(\lim_{n\to\infty} -M_{t\wedge\gamma_n}) \leq \liminf_{n\to\infty} \mathbb{E}(-M_{T\wedge\gamma_n}) = 0$$
and so $m \geq 0$ is a necessary condition for the existence of an embedding. It is also sufficient: note that $m < \alpha$ from the definition of $\mu^s$, and the stopping time $T_m = \inf\{t \geq 0 : M_t = m\}$ is a.s. finite for $m \geq 0$. We can put $T = T_m + T_{\mathbf{AY}} \circ \theta_{T_m}$, where $T_{\mathbf{AY}}$ is defined through (5.3) for the measure $\mu^s$ shifted to be centered. The stopping time $T$ embeds $\mu^s$, $M_T \sim \mu^s$. Of course, we could also use some different embedding in place of $T_{\mathbf{AY}}$.

   $\mathbf{s(I)}^\circ = (\alpha, \infty)$. A completely analogous reasoning to the one above gives that $m \leq 0$ is a necessary and sufficient condition for the existence of a solution to the embedding problem for $\mu$.

3. $\mathbf{s(I)}^\circ = (\alpha, \beta)$. In this case we have a bounded martingale and we can use the dominated convergence theorem to see that $m = \mathbb{E}\, M_T = \lim \mathbb{E}\, M_{T\wedge\gamma_n} = 0$. It is of course also a sufficient condition as we can then use any standard embedding, such as the Azéma-Yor solution (see Section 5).

We point out that the simplifying assumption we made in the first place (namely $\int |s(u)|d\mu(u) < \infty$) is not necessary. One can see that in the first case above a solution always exists and in two last cases the assumption is actually a necessary condition. We refer to Cox and Hobson [23] for more details.

## 10. Links with optimal stopping theory

There are two main ways to establish a link between an embedding problem and an optimal stopping problem. Firstly, one can prove that a solution to the optimal stopping problem is given by a certain solution to the embedding problem (for some law). Secondly, one can look for such an optimal stopping problem that its solution solves a given embedding problem. The former is discussed in Section 10.1 and the latter in Section 10.2.

Closely connected with these subjects are also the works which followed and generalized the work of Perkins in studying the interplay of laws of the supremum, the infimum and the terminal law of a martingale. The main literature in this direction is discussed in Section 10.3.

In this section we only give a brief overview of the above-mentioned topics, with no claim at describing that vast domain in a complete way. In particular we do not treat optimal stopping theory itself. For a good and brief account of its basic ideas as well as some new developments, see chapter 1 in Pedersen [86].

### *10.1. The maximality principle*

Let $\phi$ be a non-negative, increasing, continuous function and $c$ a continuous, positive function. Consider the following optimal stopping problem of maximizing

$$V_\tau = \mathbb{E}\left[\phi(S_\tau) - \int_0^\tau c(B_s)ds\right], \tag{10.1}$$

over all integrable stopping times $\tau$ such that

$$\mathbb{E}\left[\int_0^\tau c(B_s)ds\right] < \infty. \tag{10.2}$$

Suppose, in the first instance, that $\phi(x) = x$ and $c(x) = c > 0$ is constant. In this formulation the problem was solved by Dubins and Schwarz [33] in an article on Doob-like inequalities. The optimal stopping time is just the Azéma-Yor embedding for a shifted (to be centered) exponential distribution (with parameter $2c$).

The methodology allowing to deal with problems like (10.1) was presented by Dubins, Shepp and Shiryaev [34] who worked with Bessel processes and $\phi(x) = x$, $c(x) = c > 0$. The setup of $\phi(x) = x$ and $c(x)$ a non-negative, continuous function, was then treated by Peskir in a series of articles ([91], [92], [94]).

**Theorem 10.1 (Peskir [91]).** *The problem (10.1), for $\phi(x) = x$ and $c(x)$ a positive, continuous function, has an optimal solution with finite payoff if and only if there exists a maximal solution $g_*$ of*

$$g'(s) = \frac{1}{2c(g(s))(s - g(s))} \tag{10.3}$$

*which stays strictly below the diagonal in $\mathbb{R}^2$, i.e. $g_*(s) < s$ for all $s \in \mathbb{R}$. The Azéma-Yor stopping time $\tau_* = \inf\{t \geq 0 : B_t \leq g_*(S_t)\}$ is then optimal and satisfies (10.2) whenever there exists a stopping time which satisfies (10.2).*

Actually the above theorem was proven in the setup of any real, time-homogeneous diffusion to which we will come back later, as it allows us to recover the solution for a general $\phi$ as a corollary. The characterization of existence of a solution to (10.1) through existence of a solution to the differential equation (10.3) is called the *maximality principle*. We note that Dubins, Shepp and Shiryaev [34] had a different characterization of $g_*$ through the condition $\frac{g_*(t)}{t} \xrightarrow[t\to\infty]{} 1$.

Let now $\phi$ be any non-negative, non-decreasing, right-continuous function such that $\phi(B_t) - ct$ is a.s. negative on some $(t_0, \infty)$ (with $t_0$ random) and keep $c$ constant. This optimal stopping problem was solved by Meilijson [75]. Define

$$H(x) = \sup_\tau \mathbb{E}\left[\phi(x + S_\tau) - c\tau\right].$$

**Theorem 10.2 (Meilijson [75]).** *Suppose that* $\mathbb{E}\sup_t\{\phi(B_t)-ct\} < \infty$. *Then $H$ is absolutely continuous and is the minimal solution of the differential equation*

$$H(x) - \frac{1}{4c}(H'(x))^2 = \phi(x) . \qquad (10.4)$$

*If $\phi$ is constant on $[x_0,\infty)$ then $H$ is the unique solution of (10.4) that equals $\phi$ on $[x_0,\infty)$. The optimal stopping time $\tau_*$ which yields $H(0)$ is the Azéma-Yor stopping time given by $\tau_* = \inf\{t \geq 0 : B_t \leq S_t - \frac{H'(S_t)}{2c}\}$.*

Let us examine in more detail the result of Meilijson in order to compare it with the result of Peskir. The Azéma-Yor stopping time is defined as $\tau_* = \inf\{t \geq 0 : B_t \leq g(S_t)\}$ with $g(x) = x - \frac{H'(x)}{2c}$. Let us determine the differential equation satisfied by $g$. Note that $H$ is by definition non-decreasing, so we have $H'(x) = \sqrt{4c(H(x) - \phi(x))}$. Differentiating once again we obtain

$$H''(x) = \frac{2c(H'(x) - \phi'(x))}{\sqrt{4c(H(x) - \phi(x))}} = \frac{2c(H'(x) - \phi'(x))}{H'(x)}.$$

Therefore

$$g'(x) = 1 - \frac{H''(x)}{2c} = \frac{\phi'(x)}{H'(x)} = \frac{\phi'(x)}{2c(x - g(x))}. \qquad (10.5)$$

We recognize immediately the equation (10.3) only there $\phi'(s) = 1$ and $c$ was a function and not a constant. This motivated our investigation of the general problem 10.1, which is given below. It is given without proof, as it requires a fair amount of new notation and more involved ideas. The interested reader is referred to [83]. The following theorem is obtained from Theorem 3.1 in Peskir [91].

**Theorem 10.3.** *The problem (10.1) for a continuously differentiable function $\phi$ and $c$ positive, with a countable number of points of discontinuity, has an optimal solution with finite payoff if and only if there exists a maximal solution $g_*$ of*

$$g'(s) = \frac{\phi'(s)}{2c(g(s))(s - g(s))} \qquad (10.6)$$

*which stays strictly below the diagonal in $\mathbb{R}^2$, i.e. $g_*(s) < s$ for all $s \in \mathbb{R}$. In this case the payoff is given by*

$$V_* = \phi(0) - 2\int_{\phi^{-1}(g_*^Y(0))}^{\phi^{-1}(0)} uc(u)du, \qquad (10.7)$$

*where $g_*^Y$ is a function given explicitly in [83]. The stopping time $\tau_* = \inf\{t \geq 0 : B_t \leq g_*(S_t)\}$ is then optimal whenever it satisfies (10.2), otherwise it is "approximately" optimal* [10].
*Furthermore if there exists a solution $\rho$ of the optimal stopping problem (10.1)*

---

[10] In the sense of Peskir [91].

then $\mathbb{P}(\tau_* \leq \rho) = 1$ and $\tau_*$ satisfies (10.2).
If there is no maximal solution of (10.6), which stays strictly below the diagonal in $\mathbb{R}^2$, then $V_* = \infty$ and there is no optimal stopping time.

We suppose the function $\phi$ is strictly increasing and with a continuous derivative. The assumptions made by Meilijson [75] are weaker. All he assumes is that $\phi$ is non-decreasing and right-continuous. In particular he also treats the case of piece-wise constant $\phi$. The reasoning in this case is based on an iteration procedure. A similar argument can be used here to describe the stopping time and the payoff. We refer again to [83] for the details.

### *10.2. The optimal Skorokhod embedding problem*

Consider now the converse problem. That is, given a centered probability measure $\mu$ describe all pairs of functions $(\phi, c)$ such that the optimal stopping time $\tau_*$ which solves (10.1) exists, and embeds $\mu$, that is $B_{\tau_*} \sim \mu$.

This may be seen as the "optimal Skorokhod embedding problem" (Peskir [92]) as we not only specify a method to obtain an embedding for $\mu$ but also construct an optimal stopping problem, of which this embedding is a solution.

Again, consider first two special cases. First, let $\phi(x) = x$. From the Theorem 10.1 we know that $\tau_*$ is the Azéma-Yor stopping time and the function $g_*$ is just the inverse of barycentre function of some measure $\mu$. The problem therefore consists in identifying the dependence between $\mu$ and $c$. Suppose for simplicity that $\mu$ has a strictly positive density $f$ and that it satisfies the $L \log L$-integrability condition.

**Theorem 10.4 (Peskir [92]).** *In the above setup there exists a unique function*

$$c(x) = \frac{1}{2} \frac{f(x)}{\overline{\mu}(x)} \tag{10.8}$$

*such that the optimal solution $\tau_*$ of (10.1) embeds $\mu$. This optimal solution is then the Azéma-Yor stopping time given by (5.3).*

Now let $c(x) = c$ be constant. Then we have

**Theorem 10.5 (Meilijson [75]).** *In the above setup, there exists a unique function $\phi$ defined through (10.4) with $H'(x) = 2c(x - \Psi_\mu^{-1}(x))$, where $\Psi_\mu$ is the barycentre function given in (5.2), such that the optimal solution $\tau_*$ of (10.1) embeds $\mu$. This optimal solution is then the Azéma-Yor stopping time given by (5.3).*

To solve the general problem we identify all the pairs $(\phi, c)$ such that the optimal stopping time $\tau_*$ which solves (10.1) embeds $\mu$ in $B$, $B_{\tau_\mu} \sim \mu$. This is actually quite simple. We know that for $\tau_* = \inf\{t \geq 0 : B_t \leq g_*(S_t)\}$ we have $B_{\tau_*} \sim \mu$ if and only if $g_*(s) = \Psi_\mu^{-1}(s)$. This is a direct consequence of the Azéma-Yor embedding (cf. Section 5). Let us investigate the differential equation satisfied by $\Psi_\mu^{-1}$. For simplicity assume that $\mu$ has a strictly positive density $f$. It is easy to see that

$$\left(\Psi_\mu^{-1}(s)\right)' = \frac{\overline{\mu}(\Psi_\mu^{-1}(s))}{f(\Psi_\mu^{-1}(s))(s - \Psi_\mu^{-1}(s))}. \tag{10.9}$$

Comparing this with the differential equation for $g_*$ (10.6) we see that we need to have

$$\frac{\phi'(s)}{2c(\Psi_\mu^{-1}(s))} = \frac{\overline{\mu}(\Psi_\mu^{-1}(s))}{f(\Psi_\mu^{-1}(s))}. \tag{10.10}$$

Equivalently, we have

$$\frac{\phi'(\Psi_\mu(u))}{2c(u)} = \frac{\overline{\mu}(u)}{f(u)}, \quad \text{for } u \in supp(\mu). \tag{10.11}$$

If we take $\phi(s) = s$, so that $\phi'(s) = 1$, we obtain $c(u) = \frac{f(u)}{2\overline{\mu}(u)}$, a half of the so-called hazard function, which is the result of Theorem 10.4. In general we have the following theorem

**Theorem 10.6.** *Let $\mu$ be a probability measure on $\mathbb{R}$ with a strictly positive density $f$, which satisfies the $L\log L$-integrability condition. Then the optimal stopping time $\tau_*$ which solves (10.1), embeds $\mu$ in $B$, i.e. $B_{\tau_*} \sim \mu$, if and only if $(\phi, c)$ satisfies (10.11) for all $u \geq 0$. The stopping time $\tau_*$ is then the Azéma-Yor stopping time given by (5.3).*

### 10.3. Investigating the law of the maximum

Several authors investigated the laws of the maximum and the minimum of stopped processes and their dependence on the terminal and initial laws. We will not discuss here all the results as there are too many of them and they vary considerably in methods and formulations. We give instead, hopefully a quite complete, list of references together with a short description, and we leave further research to the reader.

- **Dubins - Schwarz** [33] - *this article focuses more on Doob-like inequalities but some useful descriptions and calculations for Brownian motion stopping times are provided. Formulae for constant gap stopping times are established (see Section 10.1).*

- **Kertz and Rösler** [61] - *these authors build a martingale with assigned terminal and supremum distributions. More precisely, consider a probability measure $\mu$ on $\mathbb{R}$ with $\int |x|d\mu(x) < \infty$ and any other measure $\nu$ with $\mu < \nu < \mu^*$, where $\mu^*$ is the image of $\mu$ through the Hardy-Littlewood maximal function given by (5.2). Authors show that there exists a martingale $(X_t)_{0 \leq t \leq 1}$ such that $X_1 \sim \mu$ and $S_1 \sim \nu$. It's worth mentioning that the converse problem, to prove that if $\nu \sim S_1$ then $\mu < \nu < \mu^*$, was solved much earlier in Dubins and Gilat [32] (see also Section 5.5).*

- **Jacka** [59] - *this work is concerned with Doob-like inequalities. The case of reflected Brownian motion and its supremum process is treated. Actually, a kind of maximality principle (see Section 10.1) in this special case is established.*

- **Rogers** [102] - *a characterization of possible joint laws of the maximum and the terminal value of a uniformly integrable martingale (or convergent continuous local martingale vanishing at zero) is developed. A number of previous results from the previous three articles on this list can be deduced.*

- **Vallois** [119, 120] - *the author characterized the set of possible laws of the maximum of a continuous, uniformly integrable martingale. As a corollary he obtained also a characterization of the set of possible laws of the (terminal value of) local time at* 0. *Vallois [120] also developed a solution to a particular version of the Skorokhod embedding problem, which we described in Section 3.16. In [119], independently of Rogers [102], Vallois also described the set of possible joint laws of the maximum and the terminal value of a continuous, uniformly integrable martingale.*

- **Hobson** [54] - *for fixed initial and terminal laws $\mu_0$ and $\mu_1$ respectively, the author shows that there exists an upper bound, with respect to stochastic ordering, on the law of $S_1$, where $S$ is the supremum process of a martingale $(X_t)_{0\leq t\leq 1}$ with $X_0 \sim \mu_0$ and $X_1 \sim \mu_1$. The bound is calculated explicitely and the martingale such that the bound is attained is constructed. The methods used are closely connected with those of Chacon and Walsh (see Section 3.8) and therefore with our presentation of the Azéma-Yor solution (see Section 5).*

- **Pedersen** [86] - *this Ph.D. thesis contains several chapters dealing with optional stopping and the Skorokhod problem. In particular, Doob's maximal $L^p$-inequalities for Bessel processes are discussed.*

- **Peskir** [91], [92] - *these are the works discussed above in Sections 10.1 and 10.2.*

- **Brown, Hobson and Rogers** [13] - *this paper is very similar in spirit to [54]. The authors investigate an upper bound on the law of $\sup_{0\leq t\leq 2} M_t$ with given $M_0 \equiv 0$, $M_1 \sim \mu_1$ and $M_2 \sim \mu_2$, where $(M_t)_{0\leq t\leq 2}$ is a martingale. A martingale for which this upper bound is attained is constructed using excursion theory.*

- **Hobson and Pedersen** [55] - *in this paper the authors find a lower bound on the law of the supremum of a continuous martingale with fixed initial and terminal distributions. An explicit construction, involving the Skorokhod embedding, is given. The results are applied to the robust hedging of a forward start digital option.*

- **Pedersen** [87] - *the author solves the optimal stopping problem of stopping the Brownian path as close as possible to its unknown ultimate maximum over a finite time interval.*

- **Cox and Hobson** [23] - *see Section 3.20.*

## 11. Some applications

The Skorokhod embedding problem is remarkable as it exists in the literature for over 40 years by now, but it still comes back in various forms, it poses new challenges and is used for new applications. A thorough account of applications of this problem would require at least as much space as we have taken so far to describe the problem itself and its solutions. With some regret, this is a project we shall postpone for the moment[11] . What we will do instead is try to sketch the main historical applications and give the the relevant references.

### *11.1. Embedding processes*

So far we have considered the problem of embedding, by means of a stopping time, a given distribution in a given process. Naturally, once this had been done, various applications were being considered. Primarily a question of embedding a whole process into, say, Brownian motion, was of great interest. This led to the paper of Monroe [80], where he showed that all semimartingales are in fact time-changed Brownian motions. We will not present here a full account of the domain but rather try to sketch the important steps.

A very basic question consists in representing (in law) a martingale $(X_n : n \in \mathbb{N})$ as $(B_{T_n} : n \in \mathbb{N})$, where $(B_t : t \geq 0)$ is a Brownian motion and $(T_n : n \in \mathbb{N})$ is an increasing family of stopping times. This question goes back to Harris [52] and it actually motivated both Skorokhod [115] and Dubins [31]. The question is in fact solved by any solution to the Skorokhod embedding problem. Indeed, let $(X_n : n \in \mathbb{N})$ be a martingale, $X_0 = 0$ a.s, and $T_\mu$ be some solution to the Skorokhod embedding problem. The idea is then to build consecutive stopping times $T_n$ to embed in the new starting Brownian motion $(B_{T_{n-1}+s} - B_{T_{n-1}})_{s \geq 0}$ the conditional distribution of $(X_{n+1} - X_n)$, given $X_1 = B_{T_1}, \ldots X_{n-1} = B_{T_{n-1}}$.

Formal notation is easier if we use the Wiener space. Let $\omega(t)$ be the coordinate process in the space of real, continuous functions, zero at zero, equipped with the Wiener measure, so that $(\omega(t))_{t \geq 0}$ is a Brownian motion. Fix an embedding method, which associates to a measure $\mu$ a stopping time $T_\mu$ which embeds it into $\omega$: $\omega(T_\mu(\omega)) \sim \mu$ and $\omega(T_\mu(\omega) \wedge t : t \geq 0)$ is a uniformly integrable martingale. Write $\mu_0$ for the law of $X_1$ and $\mu_n(x_1, \ldots, x_n)$ for the regular conditional distribution of $(X_{n+1} - X_n)$ given $X_k = x_k, 1 \leq k \leq n$. Let $T_1 = T_{\mu_0}$ and

$$T_n(\omega) = T_{n-1}(\omega) + T_{\mu_n(b_1,\ldots,b_n)}(\omega'_n), \qquad (11.1)$$

where $b_i = \omega(T_i(\omega))$ and $\omega'_n(s) = \omega(s + T_{n-1}(\omega)) - \omega(T_{n-1}(\omega))$. It is then straightforward to see that the following proposition holds.

**Proposition 11.1 (Skorokhod [114] p.180, Dubins [31]).** *Let $(X_n : n = 0, 1, \ldots)$ be a real valued martingale with $X_0 = 0$ a.s. and define an increasing*

---
[11]Likewise, we believe it is very instructive to see explicit forms of various solutions discussed in Section 3 for some simple probability measure $\mu$. We plan to develop this idea in a separate note.

*family of stopping times* $(T_n)$ *through* (11.1). *Then* $\omega(T_n(\omega))_{n\geq}$ *has the same distribution as* $(X_n)_{n\geq 1}$. *Furthermore if* $\mathbb{E} X_n^2 = v_n < \infty$ *for all* $n \geq 1$ *then* $\mathbb{E} T_n = v_n$ *for all* $n \geq 1$.

The next step consisted in extending this result to right-continuous martingales with continuous time parameter. It was done by Monroe [78] in the same paper, where he established the properties of minimal stopping times (see Section 8). Finally, as mentioned above, in the subsequent paper Monroe [80] established that any semimartingale is a time-changed Brownian motion. Monroe used randomized stopping times and they were obtained as a limit, so that no explicit formulae were given. Zaremba [125] tried to prove Monroe's results using only natural stopping times, but his proof is known to be false. Some years later Scott and Huggins [111], developed an embedding for reverse martingales and established for them, as an application, a functional law of iterated logarithm.

We mention two explicit embedding of processes in Brownian motion by Khoshnevisan [62] and Madan and Yor [71]. The former is an embedding of a compound Poisson process. The basic idea is as follows. First, we realize jumps independently as a sequence of i.i.d. random variables $(V_i)_{i=1}^{\infty}$. Then we define stopping times $T_i$ as consecutive excursions, of height $V_i$, of Brownian motion $B$ away from supremum, that is $T_0 = 0$, $T_{i+1} = \inf\{t > T_i : \sup_{T_i \leq u \leq t} B_u - B_t = V_{i+1}\}$. Compound Poisson process is obtained via time-changing the process $B$, so that it is linear, with given slope, in-between the jump times. Khoshnevisan imposes some conditions on the law of jumps. He then obtains Strassen's rate of convergence (as in (11.3) below), for renormalized compound Poisson process to Brownian motion, only in this case uniformly in time on compact intervals. Furthermore, he is also able to prove convergence of the number of crossings to Brownian local times[12], uniformly in space and in time (on compact intervals).

Another explicit embedding of processes in Brownian motion can be found in Madan and Yor [71]. The authors provide methods of constructing martingales with prescribed marginals. One of their techniques consists in using the Azéma-Yor solution. More precisely, if $(\mu_t)_{t\geq 0}$ is a family of zero-mean marginals they define their martingale through $M_t = B_{T_t}$, with $T_t = \inf\{u \geq 0 : S_u \geq \Psi_{\mu_t}(B_u)\}$ where $\Psi_{\mu_t}$ is the barycentre function given by (5.2). Naturally this gives the desired solution only if the family of functions $\Psi_{\mu_t}$ is increasing in $t$. The authors note that this property is equivalent to the statement that the measures $\mu_t$ are increasing in time in the mean residual life order (Shaked and Shantikumar [112], p.43).

We have considered so far embeddings in law. However, sometimes stronger a.s. representations were needed for applications. We mention a series of works by Knight ([65, 66, 67]). In his first paper, which originated from his Ph.D. dissertation, Knight [65] developed a construction of Brownian motion as a limit (uniform in finite time intervals a.s.) of random walks. He also mentioned the converse problem, i.e. the Skorokhod embedding problem, but without using this name. Knight had a good control over the rate of convergence for his limit procedure and the construction generalizes easily to arbitrary diffusions. Knight [66] then used this construction to establish the celebrated Ray-Knight theorem about the local time as Markov process in space parameter. More recently, Knight [67] used again the construction to describe the law of local time of a stopped Brownian motion, given subsequent random walk approximations. This study allowed Knight [67] to investigate the law of Brownian motion given its local time.

---

[12]The rate of convergence is seen to be $n^{1/4}(\log n)^{3/4-\epsilon}$, for any $\epsilon > 0$.

### 11.2. Invariance principles

In Skorokhod's original work [115] the embedding problem was introduced as a tool to prove some invariance principles (and their approximation rates) for random walks. This was continued afterwards and extended.

We start with a simple reasoning (Sawyer [110], Meilijson [74]) which shows how to obtain the central limit theorem for i.i.d. variables from the weak law of large numbers. Let $\mu$ be a centered probability measure with finite second moment, $\int_{\mathbb{R}} x^2 d\mu(x) = v < \infty$, and let $T$ be some solution that embeds $\mu$ in $B$, $\mathbb{E} T = v$. Write $T_1 = T$ and let $T_2$ be the stopping time $T$ for the new Brownian motion $\beta_t = B_{t+T_1} - B_{T_1}$. Apply this again to obtain $T_3$ etc. We have thus a sequence $T_1, T_2, \ldots$ of i.i.d. random variables, and $S_n = B_{\sum_{i=1}^n T_i}$ is a sum of $n$ i.i.d. variables with law $\mu$. Define the process $B_t^{(n)} = \frac{1}{\sqrt{n}} B_{nt}$, which is again a Brownian motion, thanks to the scaling property. The pairs $(B, B^{(n)})$ converge in law to a pair of two independent Brownian motions $(B, \beta)$. The stopping times $T_i$ are measurable with respect to the Brownian motion $B$ and, by the strong law of large numbers, $\frac{1}{n} \sum_{i=1}^n T_i \xrightarrow[n \to \infty]{} v$ a.s. This implies that $B^{(n)}_{\frac{1}{n} \sum_{i=1}^n T_i}$ converges in law to $\beta_v \sim \mathcal{N}(0, v)$, which ends the argument as $B^{(n)}_{\frac{1}{n} \sum_{i=1}^n T_i} = \frac{S_n}{\sqrt{n}}$. This reasoning was generalized into:

**Proposition 11.2 (Holewijn and Meilijson [56]).** *Let $X = (X_i : i \geq 1)$ be a stationary and ergodic process such that $\mathbb{E} X_1 = 0$ and $\mathbb{E}(X_n | X_1, \ldots, X_{n-1}) = 0$ a.s., $n = 2, 3, \ldots$ and $(B_t : t \geq 0)$ a Brownian motion. Then there exists a sequence of standard (randomized) stopping times $0 \leq T_1 \leq T_2 \leq \ldots$, such that*

- $(B_{T_1}, B_{T_2}, \ldots)$ *is distributed as* $X$,
- *The process of pairs* $((B_{T_1}, T_1), (B_{T_2} - B_{T_1}, T_2 - T_1), \ldots)$ *is stationary and ergodic,*
- $\mathbb{E} T_1 = Var(X_1)$.

Invariance principle is a generic name for a limit behavior of functionals of sums of random variables which proves independent of the distributions of the variables. The most known cases include the Central Limit Theorem, the Law of Iterated Logarithm and Donsker's theorem. All of these can be obtained via the Skorokhod embedding techniques, combined with appropriate properties of Brownian motion. We presented a simple example above. An advantage of this approach is that one obtains the desired limit behavior together with bounds on the rate of approximation. A survey on the subject was written by Sawyer [110] just before the Komlós-Major-Tusnády [68] presented their construction, to which we come back below.

The general setup for invariance principles can be summarized in the following way: construct a sequence of independent random vectors $X_n$, with given distributions, and a sequence of independent Gaussian vectors $Y_n$, each $Y_n$ having the same covariance structure as $X_n$, so as to minimize

$$\Delta(X,Y) := \max_{1\leq k\leq n} \Big| \sum_{i=1}^{k} X_i - \sum_{i=1}^{k} Y_i \Big|. \tag{11.2}$$

The rate of approximation for the invariance principle expresses how small $\Delta$ is as $n \to \infty$.

The classical case deals with real i.i.d. variables. Assume $X_n$ are i.i.d., $\mathbb{E}\,X_1 = 0$, $\mathbb{E}\,X_1^2 = 1$ and $\mathbb{E}\,X_1^4 = b < \infty$. To use the Skorokhod embedding to evaluate $\Delta(X,Y)$ it suffices to apply Proposition 11.1 to the martingale $M_k = \sum_{i=1}^{k} X_k$. Since the martingale increments are i.i.d. so are the consecutive stopping times' differences: $(T_{k+1} - T_k)$, $k \geq 0$, where $T_0 \equiv 0$. The Gaussian variables are given in a natural way by $Y_k = B_k - B_{k-1}$, $k \geq 1$. The differences in (11.2) reduce to $|B_{T_k} - B_k|$. In his paper Strassen [116], proved that

$$Z = \limsup_{n\to\infty} \frac{|B_{T_n} - B_n|}{(n\log\log n)^{1/4}(\log n)^{1/2}}, \tag{11.3}$$

is of order $O(1)$. Strassen used the original Skorokhod's randomized embedding, but the result does not depend on this particular construction of stopping times. This subject was further developed by Kiefer [63], who proved that the rate $((\log n)^2 n \log\log n)^{1/4}$ could not be improved under any additional moment conditions on the distribution of $X_1$. The case of independent, but not i.i.d. variables was also treated.

The Skorokhod embedding method was regarded as the best one in proving the rates of approximation in invariance principles for some time. However, in the middle of the seventies, a new method mentioned above, known as the Komlós-Major-Tusnády construction, was developed (see [25] and [68], [35] gives some extensions) that proved far better than the Skorokhod embedding and provided an exhaustive solution to minimizing (11.2). The reason why Skorokhod type embedding could not be optimal is quite clear. It was created as a means to embed a single distribution in Brownian motion and to embed a sequence one just iterates the procedure, while the Komlós-Major-Tusnády construction was created from the beginning as an embedding of random walk in Brownian motion.

Still, the Skorokhod embedding techniques can be of some use in this field. One possible example involves extending the results known for Brownian motion (as the invariance principles with $(X_n)$ uniform empirical process and $(Y_n)$ Brownian bridges) to the setup of arbitrary continuous time martingales (see for example Haeusler and Mason [48]). Also, Skorokhod embedding might be used along with the Komlós-Major-Tusnády construction (see Zhang and Hsu [127]).

### *11.3. A list of other applications*

Finally, we just point out several more fields, where Skorokhod embedding is used for various purposes. The reader is referred to the original articles for all details.

- *Financial mathematics*: some examples of applications include moving between discrete time and continuous time models (Guo [47]) and classical Black-Scholes model and other schemes (Walsh and Walsh [122]).

- *Processes in random environment:* random environment version of the Skorokhod embedding serves to embed a random walk in a random environment into a diffusion process in the random environment (see Hu [57]).

- *Study of $\frac{1}{2}$-stable distribution*: see Donati-Martin, Ghomrasni and Yor [30].

- *Invariance principle for the intersection local time* - Cadre [14] uses Azéma-Yor solution (see Section 5) to embed a two-dimensional random walk into a two-dimensional Brownian motion. He employs this embedding to prove that the renormalized intersection local time for planar Brownian motion can be obtained as a limit of the corresponding quantity in the discrete-time setting.

- *Functional CLT:* Courbot [20] uses Skorokhod embedding to study the Prohorov distance between the distributions of Brownian motion and a square integrable cádlág martingale.

**Acknowledgment.** The author would like to thank Marc Yor for introducing him to the domain and for his constant help and support which were essential from the beginning to the end of this project. The author is also thankful to Isaac Meilijson, Goran Peskir, Jim Pitman and an anonymous referee who helped greatly to improve earlier versions of this survey.